\documentclass[hidelinks,onefignum,onetabnum]{siamart220329}

\usepackage{lipsum}
\usepackage{amsfonts}
\usepackage{graphicx}
\usepackage{epstopdf}
\usepackage{algorithmic}
\usepackage{caption}
\usepackage{subcaption}
\ifpdf
  \DeclareGraphicsExtensions{.eps,.pdf,.png,.jpg}
\else
  \DeclareGraphicsExtensions{.eps}
\fi

\newsiamremark{remark}{Remark}
\newsiamremark{hypothesis}{Hypothesis}
\crefname{hypothesis}{Hypothesis}{Hypotheses}
\newsiamthm{claim}{Claim}

\headers{Fast Policy Learning for LQC with Entropy Regularization}{X. Guo, X. Li, and R. Xu}

\title{Fast Policy Learning for Linear-Quadratic Control with Entropy Regularization}

\author{Xin Guo\thanks{Department of Industrial Engineering \& Operations Research, University of California, Berkeley, Berkeley, CA 
  (\email{xinguo@berkeley.edu}).}
\and Xinyu Li\thanks{Department of Mathematics, University of Oxford, Oxford, UK
(\email{xinyu.li@maths.ox.ac.uk}).}
\and  Renyuan Xu\thanks{Department of Management Science \& Engineering, Stanford University, Stanford, CA
(\email{renyuanxu@stanford.edu}).}}

\usepackage{amsopn}

\usepackage{algorithm}
\usepackage{algorithmic}

\usepackage{mathtools}

\numberwithin{equation}{section}
\numberwithin{theorem}{section}

\newtheorem{assumption}{Assumption}[section]
\def\vec{{\text{vec}}}

\newcommand{\Tr}{\operatorname{Tr}}

\usepackage{comment}
\ifpdf
\hypersetup{
  pdftitle={Fast Policy Learning for Linear-Quadratic Control with Entropy Regularization},
  pdfauthor={Xin Guo, Xinyu Li, and Renyuan Xu }
}
\fi

\begin{document}

\maketitle

\begin{abstract}  
This paper proposes and analyzes two new policy learning methods: regularized policy gradient (RPG) and iterative policy optimization (IPO), for a class of discounted linear-quadratic control (LQC) problems over an infinite time horizon with entropy regularization. Assuming access to the exact policy evaluation,  both proposed approaches are 
proved to converge linearly in finding optimal policies of the regularized LQC.
Moreover, the IPO method can achieve a super-linear convergence rate once it enters a local region around the optimal policy. 
Finally, when the optimal policy for an RL problem with a known environment is appropriately transferred as the initial policy to an RL problem with an unknown environment,  the IPO method is shown to converge at a super-linear rate if the two environments are sufficiently close. A model-free version of the policy-based methods is also discussed.
 Performances of these proposed algorithms are supported by numerical examples.
\end{abstract}


\section{Introduction} 
Reinforcement Learning (RL) is a powerful framework for solving sequential decision-making problems, where a learning agent interacts with an unknown environment to improve her performance through trial and error \cite{DPP3_sutton2018reinforcement}.
In RL, an agent takes an action and receives
a reinforcement signal in terms of a reward, which encodes the outcome of her action. In order to maximize the accumulated reward over time, the agent learns to select her actions based on her past experiences (exploitation) and/or by making new choices (exploration).
 Exploration and exploitation are the essence of RL, and entropy regularization has shown to be effective  to 
 balance the exploration-exploitation in RL, and more importantly to enable fast convergence \cite{cen2022fast,haarnoja2018soft, mei2020global}. 
 
 Fast convergence and sample efficiency are critical for many applied RL problems, such as financial trading \cite{hambly2023recent} and healthcare treatment recommendations \cite{yu2021reinforcement}, where acquiring new samples is costly or the chance of exploring new actions in the system is limited.  In such cases, the cost of making incorrect decisions can be prohibitively high.

\paragraph{Our work}
    This paper proposes and analyzes 
two new policy learning methods: regularized policy gradient (RPG) and  iterative policy optimization (IPO),  for a  class of discounted entropy-regularized linear-quadratic control (LQC) problems over an infinite time horizon. Assuming access to the exact policy evaluation,  both approaches are shown to converge linearly in finding optimal policies of the regularized LQC  (Theorem \ref{thm: global conv of regularized PG} and \ref{thm: PI}).
Moreover, the IPO method can achieve a super-linear convergence rate (on the order of one and a half) once it enters a local region around the optimal policy. 
Finally, when the optimal policy from an RL problem with a known environment  is appropriately transferred as the initial policy to an RL problem with an unknown environment,  the IPO algorithm is shown to enable a super-linear convergence rate if the two environments are  sufficiently close
(Theorem \ref{thm:transfer_learning}).

Our analysis approach is inspired by \cite{fazel2018global} to establish the gradient dominance condition within the linear-quadratic structure. Unlike theirs, our framework incorporates entropy regularization and state transition noise (Section \ref{sec: Problem Formulation}). Therefore, 
in contrast to their deterministic and linear policies, our policies are of  Gaussian type. Consequently, the gradient dominance condition involves both the gradient of the mean and the gradient of the covariance (Lemma \ref{lemma:grad_dom}).
Accordingly, to establish the convergence of the covariance update in RPG,  the smoothness of the objective function for bounded covariance is exploited, which is ensured with proper learning rate (Lemma \ref{lemma: bound of update Sigma}).

Different from the  first-order gradient descent update in most existing literature, our proposed IPO method requires solving an optimization problem at each step. This yields faster  (super-linear) local convergence,  established by bounding the differences between two discounted state correlation matrices with respect to the change in policy parameters (Lemma \ref{lemma: local contraction of C} and Theorem \ref{thm: local conv}).
This approach is connected intriguingly  with \cite{cen2022fast}, where the bound for the difference between discounted state visitation measures yielded the local quadratic convergence in Markov Decision Processes (MDPs).

 \paragraph{Related works of policy gradient methods in LQC} As a cornerstone in optimal control theory, the LQC problem is to find an optimal control in a  linear dynamical system with a quadratic cost. LQC is popular due to its analytical tractability and its approximation power to nonlinear problems \cite{basei2022logarithmic}. 
Until recently, most works on the LQC problem assumed that the model parameters are {\it fully known}. The first global convergence
result for the policy gradient method to learn the optimal policy for LQC problems was
developed in \cite{fazel2018global} for an infinite time horizon and with deterministic dynamics. Their
work was extended in \cite{bhandari2019global}
to give global optimality guarantees of policy gradient
methods for a larger class of control problems, which satisfy a closure condition under policy
improvement and convexity of policy improvement steps. More progress has been made for policy gradient methods in other settings as well, including \cite{bu2019lqr} for
 a real-valued matrix
function,  \cite{bu2020policy} for a continuous-time setting,   \cite{gravell2020learning} for multiplicative noise, \cite{jin2020analysis}, \cite{malik2019derivative} for additive noise, and \cite{hambly2021policy} for finite time horizon with an additive noise, 
\cite{tsiamis2020risk}  and \cite{zhao2023global} for time-average costs with risk constraints, \cite{han2023policy} for nearly-linear dynamic system, 
\cite{wang2023policy} for distributional LQC to find the distribution of the return, and \cite{hamdouche2022policy} for nonlinear stochastic control with
exit time.
Our work  establishes  {\it fast} convergence  for both policy gradient based and policy optimization based algorithms for an infinite time horizon  LQC   with  entropy regularization.

\paragraph{Related works of entropy regularization}
Entropy regularization has been frequently adopted to encourage exploration and improve convergence \cite{guo2023markov,hazan2019provably,jia2022policy, jia2023q,peters2010relative,shani2020adaptive, wang2020reinforcement,wang2018exploration, williams1991function, tang2022exploratory,szpruch2022optimal,szpruch2021exploration, zhang2021policy}. In particular, \cite{ahmed2019understanding} showed that entropy
regularization induces a smoother landscape that allows for the use of larger learning rates, and hence, faster
convergence. 
Convergence rate analysis has been established when the underlying dynamic is an MDP with finite states and finite actions. For instance,   
\cite{agarwal2020optimality} and \cite{mei2020global} developed convergence guarantees for regularized policy gradient methods, with relative entropy regularization considered in \cite{agarwal2020optimality} and entropy regularization in
\cite{mei2020global}. Both papers suggest the role of regularization in guaranteeing
{\it faster} convergence for the tabular setting. 
For the natural policy gradient method, \cite{cen2022fast} established a global linear convergence rate and a local quadratic convergence rate. 

For system with infinite number of states and actions, the closest to our work in terms of model setup is \cite{wang2018exploration}. Our work replaces their aggregated control setup with controls that are randomly sampled from the policy, which are more realistic in handling real-world systems. 
The focuses of these two papers are also different: theirs explained the exploitation–exploration trade-off with entropy regularization from a continuous-time stochastic control perspective and provided theoretical support for Gaussian exploration for LQC; while ours is on algorithms design and the convergence analysis.  To the best of our knowledge,  our work is the first non-asymptotic convergence result for LQC under entropy regularization.

\paragraph{Related works of transfer learning}  
Transfer learning, \textit{a.k.a.} knowledge transfer, is a technique to utilize external expertise from other domains to benefit the learning process of a new task \cite{cao2023feasibility, cao2023risk, niu2020decade,weiss2016survey}. It has gained popularity in many areas to improve the efficiency of learning. However,  transfer learning in the RL framework is decisively more complicated and remains largely unexplored, as the knowledge to transfer involves a controlled stochastic process \cite{zhu2023transfer}. 
The transfer learning scheme proposed here is the first known theoretical development of transfer learning in the context of RL. 
% It may shed light on the benefit of integrating entropy regularization in a more general framework for transfer learning in RL.

\paragraph{Notations and organization} Throughout the paper, we will denote, for any matrix $Z \in \mathbb{R}^{m \times d}$, $Z^\top$ for the transpose of $Z$, $\|Z\|$ for the spectral norm of $Z$, $\|Z\|_F$ for the Frobenius norm of $Z$,  $\Tr(Z)$ for the trace of a square matrix $Z$, and $\sigma_{\min }(Z)$ (\textit{resp.} $\sigma_{\max}(Z)$) for the minimum (\textit{resp.} maximum) singular value of a square matrix $Z$. 
{\color{black} Let $\mathbb{S}^d_{+}$ denote the set of symmetric positive semi-definite matrices in $\mathbb{R}^{d \times d}$ and  $\mathbb{S}^d_{++}$ for the subset of $\mathbb{S}^d_{+}$ consisting of symmetric positive definite matrices.}
We will adopt  $\mathcal{N}(\mu, \Sigma)$ for a Gaussian distribution with mean $\mu \in \mathbb{R}^d$ and covariance matrix $\Sigma \in \mathbb{S}^{d}_+$.

The rest of the paper is organized as follows. Section \ref{sec: Problem Formulation} introduces the problem  and provides its theoretical solution using the dynamic programming principle.
Section \ref{sec: preliminary analysis} presents the gradient dominance condition and related smoothness property. Section \ref{sec:regularized_policy} introduces the RPG method and provides the global linear convergence result, and Section \ref{sec: policy_iteration} proposes the IPO method, along with its global linear convergence and local super-linear convergence results.  Section \ref{sec: transferlearning}  shows that IPO leads to an efficient transfer learning scheme. {A model-free version of the policy-based method is discussed in Section \ref{sec:model_free}.} Numerical examples are presented in Section \ref{sec: numerical experiments}. 
{Throughout this paper, proofs of main technical lemmas, unless otherwise specified, are deferred to Section \ref{sec: detailed proofs}.% Due to the page limit, proofs of some additional lemmas   can be found in the online companion at \url{https://arxiv.org/abs/2311.14168}.
} 

\section{Regularized LQC Problem and Solution}\label{sec: Problem Formulation}

\subsection{Problem Formulation}
We consider an entropy-regularized LQC problem over an infinite time horizon with a constant discounted rate. %\xinyu{To-do: $(\Omega, \mathcal{F}, \mathbb{P})$ set up}

\paragraph{Randomized policy and entropy regularization} To enable entropy regularization for exploration in the context of learning, we focus on randomized Markovian policies that are stationary.
%Namely, 
% {\color{black}
% Let \(\mathcal{X}\) denote the state space, \(\mathcal{A}\) the action space, and \(\mathcal{P}(\mathcal{A})\) the space of probability measures on \(\mathcal{A}\). Let \(\lambda\) represent the Lebesgue measure on \(\mathcal{A}\). We use Shannon entropy to quantify exploration: for any \(m \in \mathcal{P}(\mathcal{A})\), the entropy is defined as
% \[
% \mathcal{H}(m) = \int_{\mathcal{A}} \ln \left(\frac{m(\mathrm{d} a)}{\lambda(\mathrm{d} a)}\right) m(\mathrm{d} a),  \text{where } m \text{ is absolutely continuous with respect to } \lambda.
% \]
% The admissible policy set \(\Pi\) is defined such that for any \(\pi \in \Pi\) and \(x \in \mathcal{X}\), the conditional measure \(\pi(\cdot | x)\) is absolutely continuous with respect to \(\lambda\), i.e.,
% \[
% \mathrm{d} \pi(\cdot | x) = p_{\pi}(\cdot | x) \lambda(\mathrm{d} a),
% \]
% for some density function \(p_{\pi}(\cdot | x)\). Here, \(\pi\) maps a state \(x \in \mathcal{X}\) to a randomized action in \(\mathcal{A}\), and the corresponding Shannon entropy is given by
% \[
% \mathcal{H}(\pi(\cdot | x)) = -\int_{\mathcal{A}} \log p_{\pi}(u | x) p_{\pi}(u | x) \lambda(\mathrm{d} u).
% \]}
Namely, define the admissible policy set as $\Pi:=\{\pi: \mathcal{X} \rightarrow \mathcal{P}(\mathcal{A})\}$, with $\mathcal{X}$ the state space, $\mathcal{A}$ the action space, and $\mathcal{P}(\mathcal{A})$ the space of probability measures on action space $\mathcal{A}$. Here each admissible policy $\pi \in \Pi$ maps a state $x \in \mathcal{X}$ to a randomized action in $\mathcal{A}$.

For a given admissible policy $\pi \in \Pi$, the corresponding Shannon's entropy is defined as 
\cite{guo2022entropy,jia2023q}:
{\color{black}
\begin{eqnarray*}
    \mathcal{H}(\pi(\cdot \mid x)):=-\int_{\mathcal{A}} \log \pi (u \mid x) \pi(u \mid x) \mathrm{d} u.
\end{eqnarray*}}
 The Shannon entropy quantifies the information gain from exploring the unknown environment. We incorporate this entropy term in the objective function as a regularization to encourage collecting information in the unknown environment and performing exploration.

\paragraph{Objective function and dynamics}
The decision maker aims to find an optimal policy by minimizing the following objective function
% \begin{equation}\label{eqn: obj initial}
%      \min_{\pi\in \Pi} \mathbb{E}_{x \sim \mathcal{D}}[J_\pi(x)], \, J_\pi(x) = \mathbb{E}_\pi\left[\left. \sum_{t=1}^\infty \gamma^t \left(x_t^{\top} Q x_t + u_t^{\top} R u_t + \tau \log \pi (u_t|x_t) \right) \right| x_0 = x\right],
% \end{equation}
\begin{equation}\label{eqn: obj initial}
     \min_{\pi\in \Pi} \, \mathbb{E}_{x \sim \mathcal{D}}[J_\pi(x)], 
\end{equation}
with value function $J_\pi$ given by
\begin{equation}\label{eqn: value func}
     J_\pi(x) := \mathbb{E}_\pi\left.\left[ \sum_{t=0}^\infty \gamma^t \left(x_t^{\top} Q x_t + u_t^{\top} R u_t + \tau \log \pi (u_t|x_t) \right) \right| x_0 = x\right],
\end{equation}
%with value function $J_\pi$ given by
% \begin{equation}\label{eqn: value func}
% \end{equation}
%where $J_\pi$ denotes the value function and
and such that for $t = 0, 1, 2, \cdots$,
\begin{equation}\label{eqn: dynamic of x}
   x_{t+1} = A x_t + B u_t + w_t, \, x_0 \sim \mathcal{D}.
\end{equation}
Here $x_t \in \mathcal{X}:=\mathbb{R}^n$ is the state of the system and the initial state $x_0$ follows an initial distribution $\mathcal{D}$. Here $u_t \in \mathcal{A}:= \mathbb{R}^k$ is the control at time $t$ following a policy $\pi$. In addition, 
$\{ w_t\}_{t=0}^\infty$ are zero-mean independent and identically distributed (i.i.d) noises. We assume that  $\{ w_t\}_{t=0}^\infty$ have finite second moments. That is,  $\Tr(W) < \infty$ with $W := \mathbb{E}[w_t w_t^\top]$ for any $t = 0, 1, 2, \cdots$. 
 The matrices $A \in \mathbb{R}^{n \times n}, B \in \mathbb{R}^{n \times k}$ define the system's (transition) dynamics.
$Q \in {\color{black} \mathbb{S}^{n}_+} $ and $R \in {\color{black} \mathbb{S}_{++}^{k}}$ are matrices that parameterize the quadratic costs.
 $\gamma \in (0, 1)$ denotes the discount factor and $\tau > 0$ denotes the regularization parameter. 
The expectation in \eqref{eqn: obj initial} is taken with respect to the control $u_t\sim \pi(\cdot |x_t)$ and system noise $w_t$ for $t\ge 0$. %\rx{and system noise?}

\subsection{Optimal Value Function and Policy}
While the optimal solution to the LQC problem is a well-explored topic, it is worth noting that, to the best of our knowledge, no prior work has presented a solution to the entropy-regularized LQC problem in the form of \eqref{eqn: obj initial}. Additionally, in the study of \cite{wang2018exploration} for entropy-regularized LQC with continuous-time state dynamics,  they focused on the state transitions with aggregated controls. This differs from  the state transitions considered in \eqref{eqn: dynamic of x}, where the controls in the state transitions are randomly sampled from the policy $\pi$. 

\paragraph{Optimal value function} 
The optimal value function $J^*: \mathcal{X} \rightarrow \mathbb{R}$ is defined as
\begin{equation}\label{eqn: optimal value func}
    J^*(x) = \min_{\pi \in \Pi} J_\pi(x).
\end{equation}
The following theorem establishes the explicit expression for the optimal control policy and the corresponding optimal value function:  the optimal policy is characterized as a multivariate Gaussian distribution, with the mean linear in the state $x$ and a constant covariance matrix.
\begin{theorem}[Optimal value functions and optimal policy]\label{thm: optimal val func and policy}
The optimal policy $\pi^*$ to \eqref{eqn: optimal value func} is a Gaussian policy: $\pi^*(\cdot | x) = \mathcal{N}(-K^*x, \Sigma^*), \forall x\in \mathcal{X},$ where
  \begin{equation}\label{eqn: Sigma *}
 K^* = \gamma (R+\gamma B^\top  P B)^{-1} B^\top  P A, \quad
     \Sigma^* = \frac{\tau}{2} (R+\gamma B^\top  P B)^{-1},
 \end{equation} with $P$ and $q$ satisfying
   \begin{subequations}
\begin{align}          
P &=  Q + K^{* \top} R K^* + \gamma (A-BK^*)^{\top} P (A-BK^*) \label{eqn:P},\\
        q &=\frac{1}{1-\gamma}\Big(
        \Tr(\Sigma^* (R+\gamma B^{\top} P B)) -\frac{\tau}{2}\left(k+\log \left((2 \pi)^k \operatorname{det} \Sigma^* \right)\right)  + \gamma\Tr(W P)  \Big). \label{eqn:q}
\end{align}
\end{subequations}
The optimal value function $J^*$ in \eqref{eqn: optimal value func} can be expressed as $J^*(x) = x^\top P x + q$.
   \end{theorem}

Proof of Theorem \ref{thm: optimal val func and policy} relies on the following lemma, which establishes the optimal solution for the one-step reward function with entropy regularization in the reward.
\begin{lemma}\label{claim: solve regularzied QP} For any given matrix $M \in \mathbb{S}^{k}_+$ and vector $b \in \mathbb{R}^{k}$, the optimal solution $p^* \in \mathcal{P}( \mathcal{A})$ to the following optimization problem \eqref{eq:opt_prob} is a multivariate Gaussian distribution with  covariance $\frac{\tau}{2}M^{-1}$ and mean $-\frac{1}{2}M^{-1}b$: 
\begin{equation}\tag{P}\label{eq:opt_prob}
\begin{aligned}
    \min_{
         %p \in \mathcal{P}(\mathcal{A})
        \color{black} p: \mathcal{A} \mapsto [0, \infty)
         } &\mathbb{E}_{u \sim p(\cdot)}
         \left[
            u^{\top} M u + b^{\top} u + \tau \log p (u)
        \right],  \\
        \text{\rm subject to}  
 &\int_{\mathcal{A}} p(u)du=1. %, \qquad  p(u) \geq 0, \quad \forall u \in \mathcal{A}.
\end{aligned}
\end{equation}

\end{lemma}
\begin{proof} (of Theorem \ref{thm: optimal val func and policy}).
By the definition of $J^*$ in \eqref{eqn: optimal value func},
\begin{equation}\label{eqn: V^*}
    J^*(x) = \min_{\pi\in \Pi} \mathbb{E}_\pi
        \Bigl\{ x^{\top} Q x + u^{\top} R u + \tau \log(\pi(u|x))
        + \gamma  J^* (Ax + Bu + w) \Bigr\},
\end{equation}
where the expectation is taken with respect to $u\sim \pi(\cdot | x)$ and the noise term $w$, with mean $0$ and covariance $W$.
   Stipulating 
   \begin{equation}\label{eqn: V* stipulate}
       J^*(x) = x^\top P x + q
   \end{equation}
   for $P \in \mathbb{S}^{n}_{+}$ and $q\in \mathbb{R}$ and plugging into \eqref{eqn: V^*},
   we can obtain the optimal value function with dynamic programming principle \cite{DPP1_bellman1957markovian, DPP2_bertsekas1996neuro, DPP4_gu2023dynamic,DPP3_sutton2018reinforcement}:
    \begin{align}
    J^*(x)& =  x^{\top} Q x \notag + \min_{\pi} \mathbb{E}_\pi 
         \Bigl\{ u^{\top} R u +\tau \log(\pi(u|x)) \\
         & \quad + \gamma  \left[ (Ax + Bu + w)^{\top} P (Ax+Bu + w)  + q \right] \Bigr\}\notag \\
    & =  x^{\top} Q x + \gamma \operatorname{Tr}(WP) +  \gamma x^\top  A^\top P A x  + \gamma q \notag \\
     &\quad + \min_{\pi} \mathbb{E}_\pi \Bigl\{  u^{\top} (R+\gamma B^\top  P B) u+ \tau \log(\pi(u|x)) +  2\gamma u^\top  B^\top P A x   \Bigr\}.\label{eqn:value_func}
    \end{align}
Now apply Lemma \ref{claim: solve regularzied QP} to \eqref{eqn:value_func} with $M = R+\gamma B^\top  P B$ and $b = 2 \gamma B^\top P A x$, one can get the optimal policy at state $x$:
\begin{equation}
    \begin{aligned}\label{eqn: pi* tmp}
        \pi^*(\cdot | x) &= \mathcal{N}\left(- \gamma(R+\gamma B^\top  P B)^{-1} B^\top  P A x, \frac{\tau}{2} (R+\gamma B^\top  P B)^{-1}\right)=\mathcal{N}\left(-K^* x, \Sigma^*\right),
    \end{aligned}
\end{equation}  
where $K^*,\Sigma^*$ are defined in \eqref{eqn: Sigma *}.

To derive the associated optimal value function, we first calculate the negative entropy of policy $\pi^*$ at any state $x\in\mathcal{X}$:
\begin{equation}\label{eqn:neg_gua_entropy}
    \mathbb{E}_{\pi^*} [\log(\pi^*(u | x))] = \int_{\mathcal{A}} \log(\pi^*(u | x)) \pi^*(u | x)du = - \frac{1}{2}\left(k+\log \left((2 \pi)^k \operatorname{det} \Sigma^*\right)\right).
\end{equation}
Plug \eqref{eqn: pi* tmp} and \eqref{eqn:neg_gua_entropy} into \eqref{eqn:value_func} to get
\begin{align*}
    &J^*(x) =x^{\top}\left( Q + K^{* \top} R K^* + \gamma (A-BK^*)^{\top} P (A-BK^*) \right) x \\
         &\quad +\text{Tr}(\Sigma^* R) -\frac{\tau}{2}\left(k+\log \left((2 \pi)^k \operatorname{det} \Sigma^* \right)\right)
         +\gamma\left( \text{Tr}(\Sigma^* B^{\top} P B) + \text{Tr}(W P) + q \right).
\end{align*}
{\color{black} Combining this with \eqref{eqn: V* stipulate}, we obtain the Riccati equation in \eqref{eqn:P}, which according to \cite[Proposition 4.4.1]{bertsekas2012dynamic}  has a unique solution $P$. This is because \eqref{eqn:P} takes the  same form as the one for the classical LQR problem (without entropy regularization).} With the unique solution $P$, we can then define the unique $q$ as in  \eqref{eqn:q}, which finishes the proof.
\end{proof}

\section{Analysis of  Value Function and Policy Gradient}\label{sec: preliminary analysis}
In this section, we  analyze the expression of the policy gradient, the gradient dominance condition, and the smoothness property of the value function. These properties are necessary for studying the algorithms proposed in Section \ref{sec:regularized_policy} and Section \ref{sec: policy_iteration}.

Throughout the analysis, we assume that there exists $\rho \in (0, \frac{1}{\sqrt{\gamma}})$ satisfying $\|A-BK^* \|\leq \rho$ where $K^*$ is the optimal solution in Theorem \ref{thm: optimal val func and policy}. We consider the following domain $\Omega$ (\textit{i.e.,} the admissible control set) for $(K, \Sigma)$:
         $\Omega = \{ K \in \mathbb{R}^{k\times n}, \Sigma \in \mathbb{S}^{k}_{+}\}.$
For any $x_0 \in \mathcal{X}$ following the initial distribution $\mathcal{D}$, we assume that $\mathbb{E}_{x_0\sim \mathcal{D}}[x_0 x_0^\top]$ exists and 
$ \mu := \sigma_{\min}\left(\mathbb{E}_{x_0 \sim \mathcal{D}}[x_0 x_0^\top]\right) > 0.$
For any $(K,\Sigma) \in \Omega$, define $S_{K,\Sigma}$ as the discounted state correlation matrix, \textit{i.e.,} 
\begin{equation}\label{eqn: def S K Sigma}
    S_{K,\Sigma} := \mathbb{E}_{\pi_{K,\Sigma}}\left[ \sum_{i=0}^\infty \gamma^i x_i x_i^{\top}\right].
\end{equation}

According to Theorem \ref{thm: optimal val func and policy}, the optimal policy of \eqref{eqn: obj initial} is a Gaussian policy with a mean following a linear function of the state and a constant covariance matrix. In the remainder of the paper, we look for a parameterized policy of the form 
\begin{equation} \label{eqn: pi_simple}
    \pi_{K,\Sigma}(\cdot | x) := \mathcal{N}(-Kx, \Sigma),
\end{equation}
for any $x \in \mathcal{X}$. With a slight abuse of notation, we use $J_{K,\Sigma}$ to denote $J_{\pi_{K,\Sigma}}$ and denote the objective in \eqref{eqn: obj initial} as a function of $(K, \Sigma)$, given by 
\begin{equation}\label{eqn:C_K}
    C(K,\Sigma) := \mathbb{E}_{x \sim \mathcal{D}} \left[J_{K, \Sigma}(x)\right].
\end{equation}
 To analyze the dependence of $\Sigma$ in the objective function \eqref{eqn:C_K} for any fixed $K$, we also define $f_K: \mathbb{R}^{k\times k} \rightarrow \mathbb{R}$ as
\begin{equation}\label{eqn: f}
    f_K({\Sigma}) = \frac{\tau}{2(1-\gamma)}\log \operatorname{det} ({\Sigma})
        -\frac{1}{1-\gamma} \operatorname{Tr} \left({\Sigma}(R + \gamma B^{\top} P_K B)\right), \quad \forall \, {\Sigma} \succ 0.
\end{equation}

By applying the Bellman equation, we can get $J_{K, \Sigma}(x) = x^\top P_K x + q_{K,\Sigma}$  with matrix $P_K \in \mathbb{S}^{n}_{+}$ and $q_{K,\Sigma}\in \mathbb{R}$ satisfying
\begin{equation}\label{eqn:P_K}
\begin{aligned}
     P_K &= Q+ K^{\top}RK+\gamma(A-BK)^{\top} P_K (A-BK),\\
      q_{K,\Sigma} &= \frac{1}{1-\gamma}\Big(
        \Tr(\Sigma ( R + \gamma  B^{\top} P_K B)  -\frac{\tau}{2}\left(k+\log \left((2 \pi)^k \operatorname{det} \Sigma\right)\right) + \gamma\Tr(W P_K) 
    \Big).
\end{aligned}
\end{equation}
Note that \eqref{eqn:P_K} differs from equations \eqref{eqn:P} in terms of their solutions. In \eqref{eqn:P}, the values of $K^*$ and $\Sigma^*$ are explicitly defined by $K^* = \gamma (R+\gamma B^\top P B)^{-1} B^\top P A$ and $\Sigma^* = \frac{\tau}{2} (R+\gamma B^\top P B)^{-1}$. By substituting these values into \eqref{eqn:P}, one can obtain the solutions for $P$ and $q$, which define the optimal value function $J^*(x) = x^\top P x + q$.

Meanwhile, $K$ and $\Sigma$ in \eqref{eqn:P_K}   can take  any admissible policy parameter values in $\Omega$, and the resulting $P_K$ and $q_{K,\Sigma}$ are functions of these policy parameters. The value function $J_{K,\Sigma}(x)$ derived from \eqref{eqn:P_K} represents the value starting from state $x$ with policy parameters $(K,\Sigma)$, which may or may not correspond to an optimal policy.

We now provide an explicit form for the gradient of the objective function $C(K,\Sigma)$ with respect to $K$ and $\Sigma$. This explicit form will be used to show the gradient dominance condition in Lemma \ref{lemma:grad_dom} and also in analyzing the algorithms in Sections \ref{sec:regularized_policy} and \ref{sec: policy_iteration}.

\begin{lemma}[Explicit form for the policy gradient]\label{thm: PG}Take a policy in the form of \eqref{eqn: pi_simple} with parameter $(K,\Sigma) \in \Omega$, then 
the policy gradient has an explicit form:
    \begin{align}
        \nabla_K C (K, \Sigma) &= 2 E_K S_{K,\Sigma} , \,
        \nabla_\Sigma C(K, \Sigma)= \frac{1}{1-\gamma} \left(R - \frac{\tau}{2}\Sigma^{-1} + \gamma B^{\top} P_K B \right), \label{eqn:grad_K}
    \end{align}
where %$S_{K,\Sigma} := \mathbb{E}_\pi\left[ \sum_{i=0}^\infty \gamma^i x_i x_i^{\top} \right]$ and 
$E_K :=-\gamma B^{\top} P_K(A-BK) + R K$ and $S_{K,\Sigma}$ is defined in \eqref{eqn: def S K Sigma}.
\end{lemma}
\paragraph{Gradient dominance}
 To prove the global convergence of policy gradient methods, the key idea is to show the gradient dominance condition, which states that $C(K,\Sigma) - C(K^*,\Sigma^*)$ can be bounded by $\|\nabla_K C(K, \Sigma)\|_F^2$ and $\|\nabla_\Sigma C(K,\Sigma)\|_F^2$ for any $(K,\Sigma) \in \Omega$. 
This suggests that when the gradient norms are sufficiently small, the cost function of the given policy is sufficiently close to the optimal cost function. 

\begin{lemma}[Gradient dominance of $C(K,\Sigma)$]\label{lemma:grad_dom}Let $\pi^*$ be the optimal policy with parameters $K^*, \Sigma^*$.
    Suppose policy $\pi$ with parameter $(K, \Sigma)\in \Omega$ satisfying $\Sigma \preceq I$ has a finite expected cost, \textit{i.e.,}  $C(K,\Sigma) <\infty$. Then
    \begin{equation}\label{eqn: grad dom bound}
        C(K, \Sigma) - C(K^*, \Sigma^*) \leq 
        \frac{\| S_{K^*,\Sigma^*}\|  }{4\mu^2\sigma_{\min}(R)} \|\nabla_K C(K, \Sigma)\|_F^2
        + \frac{1-\gamma}{\sigma_{\min}(R)} \| \nabla_\Sigma C(K,\Sigma)\|_F^2.
    \end{equation}
%where $\|\cdot \|$ and $\|\cdot \|_F$ denote the operator norm and Frobenius norm of matrices respectively.

For a lower bound, with $E_K$ defined in Lemma \ref{thm: PG}, 
$$
C(K, \Sigma) - C(K^*, \Sigma^*) \geq \frac{\mu}{\| R+\gamma B^{\top} P_K     B\| }
        \operatorname{Tr}(E_K^\top E_K).
$$
\end{lemma} 
\paragraph{``Almost" smoothness}
Next, we will develop a smoothness property for the cost objective $C(K,\Sigma)$, which is necessary for establishing the convergence algorithms proposed in Section \ref{sec:regularized_policy} and Section \ref{sec: policy_iteration}.

A function $f: \mathbb{R}^n \rightarrow \mathbb{R}$ is considered {\it smooth} if the following condition is satisfied:
$
\left|f(x)-f(y)+\nabla f(x)^{\top}(y-x)\right| \leq \frac{m}{2}\|x-y\|^2,\quad \forall x, y \in \mathbb{R}^n,
$
with $m$ a finite constant \cite{fazel2018global, hambly2021policy}. 
In general, characterizing the smoothness of $C(K, \Sigma)$ is challenging, as it may become unbounded when the eigenvalues of $A - BK$ exceed $\frac{1}{\sqrt{\gamma}}$ or when $\sigma_{\min}(\Sigma)$ is close to $0$. Nevertheless, in Lemma \ref{lemma:almost_smooth}, we will see  that if $C(K, \Sigma)$ is ``almost" smooth, 
% meaning that when $K'$ and $\Sigma'$ are sufficiently close to $K$ and $\Sigma$, respectively,
then the difference $C(K, \Sigma) - C(K', \Sigma')$ can be bounded by 
%the sum of first and second-order terms 
the sum of linear and quadratic terms
involving $K - K'$ and $\Sigma - \Sigma'$. %\rx{linear and quadratic?}

\begin{lemma}[``Almost" smoothness of $C(K, \Sigma)$]\label{lemma:almost_smooth}Fix $0 < a \leq 1$ and define
$M_{a} =  \frac{ \tau\left(-\log\left(a\right) +  a-1 \right)}{ 2(1-\gamma) \left(a-1\right)^2}$. For any  $K, \Sigma$ and $K', \Sigma'$ satisfying $aI \preceq \Sigma \preceq I$ and $aI \preceq \Sigma' \preceq I$, 
\begin{equation*}
    \begin{aligned}
        &C(K', \Sigma') - C(K, \Sigma) =\operatorname{Tr}\left(S_{K',\Sigma'}(K'-K)^\top (R +\gamma B^\top P_K B)(K'-K)\right) \\
        & \quad + 2 \operatorname{Tr}\left(S_{K',\Sigma'} (K'-K)^\top E_K \right) + f_K(\Sigma) - f_K(\Sigma')
        \\
        & \leq  \operatorname{Tr}\left(S_{K',\Sigma'} (K'-K)^\top (R +\gamma B^\top P_K B)(K'-K)\right) + 2 \operatorname{Tr}\left(S_{K',\Sigma'} (K'-K)^\top E_K \right)\\
        & \quad + \frac{1}{1-\gamma}\Tr\left(  \left(R +\gamma B^\top P_K B - \frac{\tau}{2} \Sigma^{-1}\right) (\Sigma' -\Sigma)\right) + M_{a} \Tr\left( ( \Sigma^{-1}\Sigma' - I)^2\right),
        %\leq&\operatorname{Tr}\left(S_{K',\Sigma'}^x (K'-K)^\top (R +\gamma B^\top P_K B)(K'-K)\right) + 2 \operatorname{Tr}\left(S_{K',\Sigma'}^x (K'-K)^\top E_K \right) + ????\\
         %&- \frac{1}{1-\gamma}\operatorname{Tr}\left( (\Sigma - \Sigma') (R + \gamma B^{\top} P_{K} B)\right)\\
         %&+ \frac{\tau}{2(1-\gamma)} \operatorname{Tr}\left(\left((\Sigma' - \Sigma)\Sigma^{-1}\right)^2 \right) -\frac{\tau}{2(1-\gamma)} \operatorname{Tr}(\left(\Sigma'-\Sigma\right)\Sigma^{-1})
    \end{aligned}
\end{equation*}
where $f_K$ is defined in \eqref{eqn: f}.
\end{lemma}

\section{Regularized Policy Gradient Method}  \label{sec:regularized_policy}
In this section, we propose a new regularized policy gradient (RPG) update for the parameters $K$ and $\Sigma$:
\begin{equation*}
    \begin{aligned}
    K^{(t+1)} &= K^{(t)} - \eta_1 \nabla_K C\left(K^{(t)}, \Sigma^{(t)}\right)\cdot \left(S_{K^{(t)},\Sigma^{(t)}}\right)^{-1}, \\
    \Sigma^{(t+1)} &= \Sigma^{(t)} - \eta_2 \Sigma^{(t)} \nabla_\Sigma C\left(K^{(t)},\Sigma^{(t)}\right) \Sigma^{(t)}. 
\end{aligned}
\end{equation*}
RPG takes into account the inherent structure of the parameter space, which can accelerate convergence. By the explicit expressions of $\nabla_K C(K,\Sigma)$ and $\nabla_\Sigma C(K,\Sigma)$ in \eqref{eqn:grad_K}, the above update is equivalent to
\begin{equation}\tag{RPG}\label{alg: regularized update}
\begin{aligned}
     K^{(t+1)} &= K^{(t)} - 2 \eta_1 E_{K^{(t)}}, \\
     \Sigma^{(t+1)} &= \Sigma^{(t)}  - \frac{\eta_2}{1-\gamma} \Sigma^{(t)} \left( {R} + \gamma B^\top P_{K^{(t)}} B- \frac{\tau}{2} {\left(\Sigma^{(t)} \right)}^{-1}\right)\Sigma^{(t)}.
\end{aligned}
\end{equation}
From \eqref{alg: regularized update}, one can see that the update of parameter $K$ does not depend on the covariance matrix $\Sigma$. However, the update of $\Sigma$ does depend on $K$ through  $P_K$. 
%\xinyu{maybe we need to explain why the update of Sigma is in that form here}

\begin{remark}[Comparison to natural policy gradient] 
Assume that  the covariance matrix $\Sigma$ is parameterized as scalar multiplication of an identity matrix, \textit{i.e.,}  $\Sigma = sI$ for some $s >0$ and 
%\begin{equation}\label{eq:policy_diagnal}
    $\pi_{K, s} (x, u) = \mathcal{N}(Kx, sI).$
%\end{equation}
%For LQC without entropy regularization,  
Then the natural policy gradient follows the update \cite{kakade2001natural}:
\begin{eqnarray} \label{eqn:LQC_without_entropy}
    K' &= K -\eta G_K^{-1}
 \nabla C(K, sI) ,\,\,
 s' = s - \eta G_s^{-1} \partial_s C(K, sI), 
\end{eqnarray}
where $G_K$ and $G_s$ are the Fischer information matrices under policy $\pi_{K,s}$, \textit{i.e.,}
\begin{equation}
\begin{aligned}
     &G_K = \mathbb{E}\left[\sum_{t=0}^{\infty} \,\nabla \log \pi_{K,s}(u_t|x_t)\nabla \log \pi_{K,s}(u_t|x_t)^\top \right] \notag, \\
     &G_s = \mathbb{E}\left[\sum_{t=0}^{\infty} \,\partial_s \log  \pi_{K,s}(u_t|x_t) \partial_s \log  \pi_{K,s}(u_t|x_t)^\top \right].
\end{aligned}
    %\notag \label{eqn: NPG sigma update}
\end{equation}
%Here $\sigma>0$ is a fixed constant without optimization. 
When the covariance matrix of the Gaussian policy takes a diagonal form as in $\pi_{K, s}$, 
\eqref{eqn:LQC_without_entropy} are equivalent to
\begin{equation}    \label{eq:natural_policy_gradient_fixed_sigma2}
    K' = K -\eta \nabla C(K, sI) (S_{K, sI})^{-1}, \quad
     s' = s - \tilde{\eta} ~\partial_s  C(K, sI) s^2,
\end{equation}
for some constant $\tilde{\eta} > 0$.

Even though \eqref{alg: regularized update} is similar to 
\eqref{eq:natural_policy_gradient_fixed_sigma2}, there is no corresponding Fisher information form associated with \eqref{alg: regularized update} because of the additional step that simultaneously updates the covariance matrix $\Sigma$, which may not necessarily be  diagonal.
\end{remark}

We next show  that \eqref{alg: regularized update} achieves a linear convergence rate. The covariance matrices $\left\{\Sigma^{(t)} \right\}_{t=0}^\infty$ using \eqref{alg: regularized update} remain bounded, provided that the initial covariance matrix $\Sigma^{(0)}$ is appropriately bounded.
\begin{theorem}[Global convergence of \eqref{alg: regularized update}]\label{thm: global conv of regularized PG}Given $\tau \in (0, 2\sigma_{\min}(R)],$ take $(K^{(0)}, \Sigma^{(0)}) \in \Omega$ such that $\Sigma^{(0)} \preceq I$. Define  $M_\tau :=  \frac{\tau k}{2(1-\gamma)}\log \left(\frac{\sigma_{\min}(R)}{\pi\tau}\right)$ and 
$
    r_0 := \max\left\{\frac{2}{\tau\sigma_{\min}\left(\Sigma^{(0)}\right)}, \, \|R\|+ \gamma \frac{\|B^\top B\|  \left(C(K^{(0)},\Sigma^{(0)})-M_\tau\right)}{\mu + \frac{\gamma}{1-\gamma}\sigma_{\min}(W)} \right\}.
$
    Then for  $\eta_1, \eta_2 = \frac{1}{2 r_0} , \frac{\tau (1-\gamma)}{2 r_0^2}
    $
    and for 
    $
        N \geq \max \left\{
             \frac{\|S_{K^{*}, \Sigma^{*}} \| r_0}{\mu \sigma_{\min}(R) },  \frac{8 r_0^3}{\tau^2 \sigma_{\min}(R)} 
     \right\}\log \frac{C(K^{(0)},\Sigma^{(0)})-C(K^*,\Sigma^*)}{\varepsilon}, 
    $
   the regularized policy gradient descent \eqref{alg: regularized update} has the following performance bound:
   $$C\left(K^{(N)},\Sigma^{(N)}\right)-C(K^*,\Sigma^*)\leq \varepsilon.$$
\end{theorem}
{\color{black}
\begin{remark} Theorem \ref{thm: global conv of regularized PG} shows that in order to achieve an $\epsilon$-optimal value function, the number of iterations required is at least $ \mathcal{O} \left(\frac{1}{\tau^5} \right)$. Thus, the larger the value of  regularization $\tau$, the smaller the number of iterations, and the faster the convergence. \end{remark}

\begin{remark}
    In the RPG update, we only need an upper bound and a lower bound for $\Sigma$, namely, $a I \preceq \Sigma \preceq b I$ for some $0 < a< b.$ Different choices of $b$ may lead to different (admissible) ranges for $\tau$.  For ease of exposition, we set $b = 1$ in Theorem \ref{thm: global conv of regularized PG}, and the results can be easily extended to the general case of any $b \geq a > 0$.
\end{remark} }

To prove Theorem \ref{thm: global conv of regularized PG}, we will first need  the boundedness of the one-step update of the covariance $\Sigma$, in order to guarantee the well-definedness of the objective function along the trajectory when performing \eqref{alg: regularized update} (Lemma \ref{lemma: bound of update Sigma}). Additionally, we will bound the one-step update of \eqref{alg: regularized update} (Lemma \ref{lemma:PG contraction}), and provide an upper bound of $\|P_K\|$ in terms of the objective function $C(K,\Sigma)$ (Lemma \ref{lemma: bounded of PK}).

\begin{lemma}[Boundedness of the update of $\Sigma$ in \eqref{alg: regularized update}]\label{lemma: bound of update Sigma}Let $(K, \Sigma)\in\Omega$ be given such that $0 \prec \Sigma \preceq I$. Fix $\tau \in (0, \, 2 \sigma_{\min}(R)]$, $a \in \Big(0, \,  \min \Big\{ \frac{\tau}{2 \| R+\gamma B^\top P_{K}B\|},  \sigma_{\min}(\Sigma) \Big\}\Big).$
    Let $K', \Sigma'$ be the one-step update of $K, \Sigma$ using \eqref{alg: regularized update} with $\eta_2 \leq \frac{2(1-\gamma)a^2 }{\tau}$. Then 
    $aI \preceq \Sigma' \preceq I$.
\end{lemma}

\begin{lemma}[Contraction of \eqref{alg: regularized update}]
\label{lemma:PG contraction}Let $(K, \Sigma)\in\Omega$ be given such that $0 \prec \Sigma \preceq I$. Assume $\tau \in (0, 2 \sigma_{\min}(R)]$.
    Fix $a \in \left(0, \, \min \left\{ \frac{\tau}{2 \| R+\gamma B^\top P_{K}B\|}, \sigma_{\min}(\Sigma) \right\}\right).$
    Let $K', \Sigma'$ be the one-step update of $K, \Sigma$ using \eqref{alg: regularized update} with $\eta_1 \leq\frac{1}{2 \| R +\gamma B^\top P_K B \|}, \eta_2 \leq \frac{2(1-\gamma)a^2}{\tau}$. Then 
    $aI \preceq \Sigma' \preceq I$
and
$$C(K', \Sigma') -C(K^*, \Sigma^*)\leq (1-\zeta) (C(K, \Sigma) -  C(K^*, \Sigma^*)),$$
with 
$ 0 < \zeta = \min \left\{\frac{2\mu\eta_1 \sigma_{\min}(R)}{\| S_{K^*,\Sigma^*} \|}, \frac{\eta_2 a \sigma_{\min}(R)}{2 (1-\gamma)}  \right\} <1.$
\end{lemma}
\begin{lemma}[Lower bound for $C(K, \Sigma)$]\label{lemma: bounded of PK}Let $M_\tau$ be defined in the same way as in Theorem \ref{thm: global conv of regularized PG}.
%There exists a constant $M_\tau =  \frac{\tau k}{2(1-\gamma)}\log (\frac{\sigma_{\min}(R)}{\pi\tau})$  such that for all $K, \Sigma$, 
Then for all $(K, \Sigma) \in \Omega$,
$C(K,\Sigma)  \geq  \left(\mu + \frac{\gamma}{1-\gamma} \sigma_{\min}(W)\right) \| P_K\| + M_\tau.$
\end{lemma}
  \begin{proof} (of Theorem \ref{thm: global conv of regularized PG}).
 Using Lemma \ref{lemma: bounded of PK} for any $t \geq 0,$
 \begin{equation}\label{eqn: eta 1 in a thm}
     \frac{1}{\left\|R+\gamma B^{\top} P_{K^{(t)}} B\right\|}
 \geq \frac{1}{\|R\|+\gamma \|B^\top B\| \left\|P_{K^{(t)}}\right\|}
 \geq \frac{1}{\|R\|+\gamma \frac{\|B^\top B\| \left( C({K^{(t)}}, \Sigma^{(t)})-M_\tau \right)}{\mu + \frac{\gamma}{1-\gamma}\sigma_{\min}(W)}}.
 \end{equation}
 Let $a=\frac{\tau}{2 r_0}$, $\zeta = \min \left\{ \frac{2\mu\eta_1 \sigma_{\min}(R)}{\|S_{K^*,\Sigma^*} \|} , \frac{\eta_2 a \sigma_{\min}(R)}{2(1-\gamma)} \right\}$. The proof is completed by induction to show $C(K^{(t+1)}, \Sigma^{(t+1)}) \leq (1-\zeta)C(K^{(t)}, \Sigma^{(t)})$ and $aI \preceq \Sigma^{(t+1)} \preceq I$ holds for all $t\geq 0$: at $t=0$, apply \eqref{eqn: eta 1 in a thm} to get $\eta_1 \leq \frac{1}{2 \| R+ \gamma B^\top P_{K^{(0)} } B\|}$ and $a \leq 
 \frac{\tau}{2\| R+ \gamma B^\top P_{K^{(0)} } B\| }$. Additionally with $ \eta_2 = \frac{\tau (1-\gamma)}{2 r_0^2}=\frac{2(1-\gamma) a^2}{\tau}$ and $aI \preceq \Sigma^{(0)} \preceq I$, Lemma \ref{lemma:PG contraction} can be applied to get
 $C(K^{(1)}, \Sigma^{(1)}) \leq (1-\zeta) C(K^{(0)}, \Sigma^{(0)}) \leq C(K^{(0)}, \Sigma^{(0)}),$
 and $aI \preceq \Sigma^{(1)} \preceq I$.
 The proof proceeds by arguing that Lemma \ref{lemma:PG contraction} can be applied at every step. If it were the case that $C\left(K^{(t)},\Sigma^{(t)}\right)\leq (1-\zeta) C\left(K^{(t-1)},\Sigma^{(t-1)}\right)\leq C\left(K^{(0)}, \Sigma^{(0)}\right)$ and $aI \preceq \Sigma^{(t)} \preceq I$, then 
% \begin{align*}
$$     2 \eta_1 
 = 
    \frac{1}{\|R\|+ \gamma \frac{\|B^\top B\|   \left(C(K^{(0)},\Sigma^{(0)})-M_\tau\right)}{\mu + \frac{\gamma}{1-\gamma}\sigma_{\min}(W)}} 
\leq  
    \frac{1}{\|R\|+ \gamma \frac{\|B^\top B\| \left(C\left(K^{(t)},\Sigma^{(t)}\right)-M_\tau\right)}{\mu + \frac{\gamma}{1-\gamma}\sigma_{\min}(W)}},$$ 
% \end{align*}
    thus by \eqref{eqn: eta 1 in a thm} $\eta_1 \leq \frac{1}{\left\|R+\gamma B^{\top} P_{K^{(t)}} B\right\|}$ and in the same way  $a\leq \frac{\tau}{2 \left\|R + \gamma B^\top P_{K^{(t)}}B\right\|}$. 
    Thus, Lemma \ref{lemma:PG contraction} can be applied such that
$
    C\left(K^{(t+1)},\Sigma^{(t+1)}\right) -  C(K^*,\Sigma^*) \leq (1-\zeta) \left(C\left(K^{(t)},\Sigma^{(t)}\right)-C(K^*,\Sigma^*)\right),
$
and $aI \preceq \Sigma^{(t+1)} \preceq I.$
Thus, the induction is complete.
Finally, observe that
    $0 < \zeta \leq \frac{2\mu\eta_1 \sigma_{\min}(R)}{\|S_{K^{*}, \Sigma^{*}} \|}  = \frac{\mu \sigma_{\min}(R)}{\|S_{K^{*}, \Sigma^{*}} \| r_0}< 1, \text{ and }
    \zeta \leq \frac{\eta_2 a \sigma_{\min}(R)}{2(1-\gamma)} 
    = \frac{\tau^2 \sigma_{\min}(R)}{8 r_0^3},$
the proof is complete.
\end{proof}
\section{Iterative Policy Optimization Method}\label{sec: policy_iteration}
In this section, we propose  an iterative policy optimization method (IPO), in which
we optimize a one-step deviation from the current policy in each iteration. For IPO, one can show both the global convergence with a linear rate and  a local super-linear convergence when the initialization is close to the optimal policy. This local super-linear convergence result benefits from the entropy regularization.

\begin{comment}
From \eqref{eqn:P_K}, one can obtain the value function starting at $x$ with parameter $(K, \Sigma):$
\begin{align*}
    J_{K, \Sigma}(x) =  &x^\top (Q+K^\top R K) x + \Tr(\Sigma R) - \frac{\tau}{2} \left( k + \log\left( (2\pi)^k\operatorname{det}\Sigma \right)\right)\\
    &
    + \gamma x^\top (A- BK)^\top P_K (A - BK) x + \gamma \Tr(W P_K) + \gamma \Tr(\Sigma B^\top P_K B).
\end{align*}
We assume the one-step update $(K', \Sigma')$ satisfies:
\begin{align*}
    K',\Sigma' = \arg\min_{\widetilde{K}, \widetilde{\Sigma}} \, &x^\top (Q+\widetilde{K}^\top R \widetilde{K}) x + \Tr(\widetilde{\Sigma} R) - \frac{\tau}{2} \left( k + \log\left( (2\pi)^k\operatorname{det}\widetilde{\Sigma} \right)\right)\\
    &
    + \gamma x^\top (A- B\widetilde{K})^\top P_K (A - B\widetilde{K}) x + \gamma \Tr(W P_K) + \gamma \Tr(\widetilde{\Sigma} B^\top P_K B).
\end{align*}
\end{comment}
{ By the Bellman equation for the value function $J_{K,\Sigma}$,
\begin{align*}
        J_{K,\Sigma}(x)= \mathbb{E}_{u\sim \pi_{K,\Sigma}} \Big[ & x^\top Q  x + u^\top R u +\tau \log \pi_{K,\Sigma}(u|x) + \gamma J_{K,\Sigma}(Ax + Bu + w)\Big].
\end{align*}
We assume the one-step update $(K', \Sigma')$ satisfies:
\begin{align*}
        K',\Sigma' = \arg \min_{\widetilde{K},\widetilde{\Sigma}}  \mathbb{E}_{u\sim \pi_{\widetilde{K},\widetilde{\Sigma}}} \Big[  x^\top Q  x + u^\top R u +\tau \log \pi_{\widetilde{K},\widetilde{\Sigma}}(u|x) + \gamma { J_{K,\Sigma}}(Ax + Bu + w)\Big].
    \end{align*} }
% and 
% $$
%     \Sigma' = \arg \min_{ \widetilde\Sigma} \,\Tr(\widetilde\Sigma R) - \frac{\tau}{2} \left( k + \log\left( (2\pi)^k\operatorname{det}\widetilde\Sigma \right)\right) + \gamma \Tr(\widetilde\Sigma B^\top P_K B).
% $$
By direct calculation, we have the following explicit forms for the updates:
\begin{equation} \tag{IPO}\label{alg: PI}
    \begin{aligned}
      K^{(t+1)} &= K^{(t)} - \left(R + \gamma B^\top P_{K^{(t)}} B\right)^{-1} E_{K^{(t)}}, \\
      \Sigma^{(t+1)} &= \frac{\tau}{2} \left( R + \gamma B^\top P_{K^{(t)}} B\right)^{-1},
\end{aligned}
\end{equation}
for $t=1,2,\cdots.$ The update of $K$ in \eqref{alg: PI} is identical to the Gauss-Newton update when the learning rate is equal to $1$ in \cite{fazel2018global}. The update of $\Sigma$ in \eqref{alg: PI} is not gradient-based and only depends on the value of $K$ in the previous step.

\subsection{Global Linear Convergence}\label{sec: PI}
In this section, we establish the global convergence for \eqref{alg: PI} with a linear rate. 

\begin{theorem}[Global convergence of \eqref{alg: PI}] \label{thm: PI}For 
$$
N \geq \frac{\left\|S_{K^*, \Sigma^*} \right\|}{\mu} \log  \frac{C(K^{(0)}, \Sigma^{(0)} )-C(K^*, \Sigma^*)}{\varepsilon},
$$ the iterative policy optimization algorithm \eqref{alg: PI} has the following performance bound:
$$
C(K^{(N)}, \Sigma^{(N)}) - C(K^*, \Sigma^*) \leq \varepsilon.
$$
\end{theorem}
The proof of Theorem \ref{thm: PI} is immediate given the following lemma, which bounds the one-step progress of \eqref{alg: PI}:
\begin{lemma}[Contraction of \eqref{alg: PI}]\label{lemma: contraction PI}Suppose $(K',\Sigma')$ follows a one-step updating rule of \eqref{alg: PI} from $(K,\Sigma)$.
Then
$$
       C(K', \Sigma') - C(K^*, \Sigma^*) \leq  \left(1-\frac{\mu}{\| S_{K^*, \Sigma^*}\|}\right) \left(C(K,\Sigma) -  C(K^*,\Sigma^*)  \right),
$$
with $ 0 < \frac{\mu}{\| S_{K^*, \Sigma^*}\|} \le 1$.
\end{lemma}

%\begin{proof}{\bf of Theorem \ref{thm: PI}} The global convergence of \eqref{alg: PI} is immediate with Lemma \ref{lemma: contraction PI} \end{proof}
Theorem \ref{thm: PI} suggests that \eqref{alg: PI} achieves a global linear convergence. Compared with 
\eqref{alg: regularized update}, \eqref{alg: PI} exhibits faster convergence in terms of the rate at which the objective function decreases (\textit{cf.} Lemmas \ref{lemma:PG contraction} and \ref{lemma: contraction PI}).
Furthermore, the subsequent section demonstrates that \eqref{alg: PI} enjoys a local super-linear convergence when the initial policy parameter is within a neighborhood of the optimal policy parameter.

\subsection{Local Super-linear Convergence}\label{sec: local super-linear convergence}
This section establishes a local super-linear convergence for \eqref{alg: PI}. We first introduce some constants used throughout this section:
\begin{equation}\label{eqn: constants xi zeta omega}
    \begin{aligned}
    \xi_{\gamma, \rho} &:= \frac{1-\gamma \rho^2 +\gamma}{(1-\gamma \rho^2)^2}, \qquad
    \zeta_{\gamma, \rho} :=  \frac{2-\rho^2}{(1-\rho^2)^2(1-\gamma)}+ \frac{1}{(1-\rho^2)^2  (1-\gamma \rho^2)},\\
    \omega_{\gamma, \rho} &:=  \frac{1}{(1-\rho^2) (1-\gamma)} + \frac{1}{(1-\rho^2)(1-\gamma \rho^2)}.
\end{aligned}
\end{equation}
To simplify the exposition, we often make use of the notation $S_{K^{(t)}, \Sigma^{(t)}}$ and $S_{K^*, \Sigma^*}$, which we abbreviate as $S^{(t)}$ and $S^*$, respectively, provided that the relevant parameter values are clear from the context.
Then, define
\begin{equation}\label{eqn: define C1 C2}
    \begin{aligned}
        \kappa:=  &\frac{\rho + \| A\|}{|\sigma_{\min}(B)|}, \quad c :=  2\rho \xi_{\gamma \rho}  \|B\| (\|Q\| + \| R\|\kappa^2 ) + \frac{1}{\mu} \|S^*\| \| R \| \cdot ( \kappa  + \| K^*\| ), \\
         c_1 := &\left( \xi_{\gamma, \rho} 
     \| \mathbb{E}[x_0 x_0^\top]\|+ \zeta_{\gamma, \rho}  \|B\Sigma^* B^\top +W\| \right) \cdot 2\rho \|B\|\\
      & \cdot \Big( 1+ \sigma_{\min}(R) \cdot \|R +\gamma B^{\top} P_{K^*}B  \|  + c \gamma  \sigma_{\min}(R) \cdot \left( \| B\| \| A\| + \| B\|^2 \kappa \right)  \Big), \\
      c_2 := &     \frac{c \tau \gamma\omega_{\gamma, \rho}  \|B\|^4}{ 2  \sigma_{\min}(R)^2}.
    \end{aligned}
\end{equation}   
Note that for any $K \in \Omega$, $\| K\| \leq \frac{\| BK\|}{|\sigma_{\min}(B)|}\leq \frac{\| A- BK\| + \| A\|}{|\sigma_{\min}(B)|} \leq \frac{\rho + \| A\|}{|\sigma_{\min}(B)|} = \kappa.$  

 We now show that \eqref{alg: PI} achieves a super-linear convergence rate once the policy parameter $(K, \Sigma)$ enters a neighborhood of the optimal policy parameter $(K^*, \Sigma^*).$

\begin{theorem}[Local super-linear convergence of \eqref{alg: PI}]\label{thm: local conv}Let $c_1$ and $c_2$ be defined in \eqref{eqn: define C1 C2}. Let $\delta := \min \left\{ \frac{1}{c_1 + c_2} \sigma_{\min}(S^*) , \frac{\rho - \|A -BK^* \|}{\| B\|} \right\}$.
Suppose that the initial policy $(K^{(0)}, \Sigma^{(0)})$ satisfies
\begin{equation}\label{eqn: assumption on C K  - C *}
    C(K^{(0)}, \Sigma^{(0)}) - C(K^*,\Sigma^*) \leq \left( \frac{1}{\mu} - \frac{1}{\|S^*\|}\right)^{-1}\sigma_{\min}\left( R + \gamma B P_{K^*} B\right) \delta^2,
\end{equation}
then the iterative policy optimization algorithm \eqref{alg: PI} has the following convergence rate: for $t = 0,1,2, \cdots$,
\begin{equation*}
     C(K^{(t+1)},\Sigma^{(t+1)}) - C(K^*, \Sigma^{*})
        \leq  \frac{ (c_1 +  c_2) \left(C(K^{(t)}, \Sigma^{(t)}) - C(K^*, \Sigma^{*}) \right)^{1.5}}{\sigma_{\min}(S^*) \sqrt{\mu \sigma_{\min}(R+\gamma B^\top P_{K^*}B)}}.
\end{equation*}   
\end{theorem}
The following Lemma \ref{lemma: local contraction of C} is critical for establishing this local super-linear convergence: it shows that there is a contraction if the differences between two discounted state correlation matrices $ \|S^{(t+1)} - S^*\|$ is small enough. 
Then, by the perturbation analysis   for $S_{K, \Sigma}$ (Lemma \ref{lemma: Sx perturbation}), one can bound 
$\|S^{(t+1)} - S^*\|$ by $\| K^{(t)} - K^*\|$ (Lemma \ref{lemma: Sprime minus Sstar}). The proof of Theorem \ref{thm: local conv} follows by ensuring the admissibility of model parameters $\left\{K^{(t)}, \Sigma^{(t)}\right\}_{t=0}^\infty$ along all the updates.

% Moreover, it is necessary to establish the following   perturbation analysis for $S_{K, \Sigma}$ (Lemma \ref{lemma: Sx perturbation}), and the bound of   $\|S^{(t+1)} - S^*\|$ by $\| K^{(t)} - K^*\|$ (Lemma \ref{lemma: Sprime minus Sstar}).
\begin{lemma} \label{lemma: local contraction of C}
Suppose that 
$\|S^{(t+1)} - S^* \|  \leq \sigma_{\min}(S^*)$ for all $t\geq 0$ when updating with \eqref{alg: PI},  then 
\begin{equation*}
    C(K^{(t+1)}, \Sigma^{(t+1)}) - C(K^*, \Sigma^*) \leq \frac{\| S^{(t+1)} - S^* \| }{\sigma_{\min}(S^*)} \left( C(K^{(t)}, \Sigma^{(t)}) - C(K^*, 
        \Sigma^*)\right).
\end{equation*}
\end{lemma}
\begin{lemma}[$S_{K, \Sigma}$ perturbation]\label{lemma: Sx perturbation}
For any $(K_1, 
\Sigma_1)$ and $(K_2, \Sigma_2)$ satisfying  $\|A - BK_1 \| \leq \rho$ and $\|A - BK_2 \| \leq \rho$, 
\begin{align*}
      & \|S_{K_1, \Sigma_1} - S_{K_2, \Sigma_2}\| \\
      \leq & \left(\xi_{\gamma, \rho}  
     \| \mathbb{E}[x_0 x_0^\top]\|+\zeta_{\gamma, \rho}  \|B\Sigma_1 B^\top +W\| \right) \cdot 2 \rho \|B\|\left\|K_1-K_2\right\| +  \omega_{\gamma, \rho}\| B\|^2 \| \Sigma_1 - \Sigma_2\|. 
\end{align*}
\end{lemma}

\begin{lemma}[Bound of one-step update of $S^{(t)}$]\label{lemma: Sprime minus Sstar}Assume that the update of parameter $(K, \Sigma)$ follows \eqref{alg: PI}.
Let $c_1, c_2$ be defined in \eqref{eqn: define C1 C2}.
Then for $K^{(t+1)}$ satisfying $\|A-BK^{(t+1)} \| \leq \rho$, we have
  $
        \|S^{(t+1)} - S^* \| \leq (c_1 + c_2) \| K^{(t)} - K^*\|.
$
\end{lemma}
\begin{proof} (of Theorem \ref{thm: local conv}).
First, Theorem \ref{thm: optimal val func and policy} shows that for an optimal $K^*$, $K^* = \gamma(R+\gamma B^{\top} P_{K^*} B)^{-1} B^{\top} P_{K^*} A.$
Then, by the definition of $E_K$ in Lemma \ref{thm: PG},
%Apply \eqref{eqn: E K} with $K^*$ to get:
\begin{equation}\label{eqn: E_K *}
    E_{K^*} =-\gamma B^{\top} P_{K^*}A+ (\gamma B^{\top} P_{K^*}B + R )K^* = 0.
\end{equation}
Fix integer $t \geq 0$.  Observe that
\begin{equation}\label{eqn: C- C* vs K- K*}
   \begin{aligned}
 & \left(1-\frac{\mu}{\| S_{K^*, \Sigma^*}\|}\right) \left(C(K^{(t)}, \Sigma^{(t)}) - C(K^*, \Sigma^*) \right) 
\overset{(a)}{\geq}
 C(K^{(t+1)}, \Sigma^{(t+1)}) - C(K^*, \Sigma^*) 
 \\
 &\overset{(b)}{\geq}   \mu \sigma_{\min}\left( R + \gamma B P_{K^*} B\right) \| K^{(t+1)}  - K^*\|^2  +f_{K^*}(\Sigma^*) - f_{K^*}(\Sigma^{(t+1)}) \\
 & \overset{(c)}{\geq} \mu \sigma_{\min}\left( R + \gamma B P_{K^*} B\right) \| K^{(t+1)} - K^*\|^2.
\end{aligned} 
\end{equation}
 $(a)$ is from  the contraction property in Lemma \ref{lemma: contraction PI};  $(b)$ follows from Lemma \ref{lemma:almost_smooth} and \eqref{eqn: E_K *};  to obtain $(c)$, note that
$
f_{K^*}\left(\Sigma^*\right) - f_{K^*}\left(\Sigma^{(t+1)}\right) \geq 0,
$
since $\Sigma^*$ is the maximizer of $f_{K^*}.$
% With % the assumption \rx{assumption?} in
Thus, \eqref{eqn: C- C* vs K- K*} and \eqref{eqn: assumption on C K  - C *}  imply
$
 \| K^{(t+1) } - K^*\| \leq \delta %\min \Big\{ \frac{1}{c_1 +  c_2} \sigma_{\min}(S^*) , \frac{\rho - \|A -BK^* \|}{\| B\|} \Big\},
$
which suggests that
$ \| A - BK^{(t+1)}\| \leq \| A - B K^*\| + \| B\| \|K^{(t+1)} - K^* \| \leq \rho.$
%This inequality along with  Lemma \ref{lemma: Sprime minus Sstar} leads 
Then by Lemma \ref{lemma: Sprime minus Sstar},
\begin{equation} \label{eqn: bound S}
    \| S^{(t+1)} - S^{*}\| \leq (c_1 + c_2) \| K^{(t)} - K^*\| \leq \sigma_{\min}(S^*).
\end{equation}
Thus, one can apply Lemma
\ref{lemma: local contraction of C} to get: 
\begin{equation}\label{eqn: C' - C local conv}
        C(K^{(t+1)}, \Sigma^{(t+1)}) - C(K^*, \Sigma^*) 
        \leq \frac{c_1 +  c_2}{\sigma_{\min}(S^*)}  \| K^{(t)} - K^*\|(C(K^{(t)},  \Sigma^{(t)}) - C(K^*, \Sigma^*)).
\end{equation}
Using the same reasoning as in \eqref{eqn: C- C* vs K- K*} $(a)$ to $(c)$, we have
$C(K^{(t)}, \Sigma^{(t)}) - C(K^*, \Sigma^*)  {\geq} \\ \mu \sigma_{\min}\left( R + \gamma B^\top P_{K^*} B\right) \| K^{(t)} - K^*\|^2$ and plug it
 in \eqref{eqn: C' - C local conv} finishes the proof.
   \end{proof}

\section{Transfer Learning for RL} \label{sec: transferlearning}One can apply the local super-linear convergence result in Theorem \ref{thm: local conv}  to provide an efficient policy transfer from a well-understood environment to a new yet similar environment. The idea is to use the optimal policy from the well-understood environment as an initialization of the policy update. If this initial policy is within the super-linear convergence region of the new environment, one may  efficiently learn the optimal policy in the new environment.  
\paragraph{Problem set-up and main results} 
We analyze two environments $\mathcal{M}:=(A,B)$  and  $\overline{\mathcal{M}}:=(\overline{A},\overline{B})$, with $(K^*,\Sigma^*)$  and $(\overline{K}^*,\overline{\Sigma}^*)$ as their respective  optimal policies and $C$ and $\overline{C}$ as their respective objective functions. Assume that one has access to the optimal (regularized) policy $(K^*,\Sigma^*)$ for environment $\mathcal{M}$, called the {\it well-understood environment}. We use $(K^*,\Sigma^*)$ as a {\it policy initialization} for the less understood environment $\overline{\mathcal{M}}$, called  the {\it new environment}. The goal is to investigate under what conditions this initialization enters the super-linear convergence regime of $\overline{\mathcal{M}}$.

Throughout this section, we specify the operator norm $\|\cdot\|$ as the one associated with vector $q$-norm. Namely, for $q\in(0,1)$ and $A\in \mathbb{R}^{n_1\times n_2}$ for some positive integers $n_1,n_2$:
$\|A\|:=\|A\|_q = \sup_{x \neq 0}\left\{\frac{\|Ax\|_q}{\|x\|_q}, \,\,  x\in \mathbb{R}^{n_2} \right\}.$
For ease of the analysis and to make the two environments comparable,  the following assumptions are made:
\begin{assumption} Assume the following conditions hold: \\
1.  Admissibility: $(K^*,\Sigma^*)$ is admissible for $\overline{\mathcal{M}}$ and $\left(\overline{K}^*,\overline{\Sigma}^*\right)$ is admissible for ${\mathcal{M}}$, \textit{i.e.,}  $\| A - BK^*\| \leq \rho, \| \overline{A} - \overline{B} \,\overline{K}^*\| \leq \rho$ with $\rho \in (0, \frac{1}{\sqrt{\gamma}})$ and $\Sigma^* \succeq 0, \overline{\Sigma}^* \succeq 0.$\\
2.  Model parameters:
 $\|B\|_q,\left\|\overline{B}\right\|_q \leq 1$.\\
3.  Optimal policy: $\|K^*\| \leq 1$.
\end{assumption}
The first condition ensures that the environments $\mathcal{M}$ and $\overline{\mathcal{M}}$ are comparable. The second and third conditions are for ease of  exposition and can be easily relaxed.

Similar to $P_K$ defined in \eqref{eqn:P_K} for environment $\mathcal{M}$, let us define the Riccati equation for  the new environment $\overline{M}$ as:
\begin{eqnarray}
  %  P_K &=& \gamma(A-BK)^{\top} P_K (A-BK) + Q+ K^{\top}RK\label{eqn:P_K1}\\
    \overline{P}_K &=& \gamma(\overline{A}-\overline{B}K)^{\top} \overline{P}_K (\overline{A}-\overline{B}K) + Q+ K^{\top}RK,\label{eqn:P_K2}
\end{eqnarray}
and define $\kappa', c', c_1', c_2'$ in the same way as \eqref{eqn: define C1 C2} with $(\overline{A}, \overline{B})$ replacing $(A, B)$.
% Let $\kappa' := \frac{\rho + \| \overline{A}\|}{|\sigma_{\min}(\overline{B})|}$ and  define the following constants:
% \begin{equation}\label{eqn: define C1 C2 prime}
%     \begin{aligned}
%         c' := & 2\rho \xi_{\gamma \rho}  \|\overline{B}\| (\|Q\| + \| R\|(\kappa')^2 ) + \frac{1}{\mu} \|S_{\overline{K}^*,\overline{\Sigma}^*}\| \| R \| \cdot ( \kappa'  + \| \overline{K}^*\| ), \\
%          c'_1 := & \left( \xi_{\gamma, \rho} 
%      \| \mathbb{E}[x_0 x_0^\top]\|+ \zeta_{\gamma, \rho}  \|\overline{B}\overline{\Sigma}^* \overline{B}^\top +W\| \right) \cdot 2\rho \|\overline{B}\|\\
%      &\cdot \Big( 1+ \sigma_{\min}(R)  \|\gamma \overline{B}^{\top} \overline{P}_{\overline{K}^*}\overline{B} + R \|  + c' \gamma  \sigma_{\min}(R) \cdot \left( \| \overline{B}\| \| \overline{A}\| + \| \overline{B}\|^2 \kappa' \right)  \Big) \\
%       c'_2 := &     \frac{c' \omega_{\gamma, \rho} \tau \gamma \|\overline{B}\|^4}{ 2  \sigma_{\min}(R)^2}.
%     \end{aligned}
% \end{equation}   

The following theorem suggests that if the environments $\mathcal{M}$ and $\overline{\mathcal{M}}$ are sufficiently close in the sense of \eqref{eq:close}, then $(K^*,\Sigma^*)$ serves an efficient initial policy for $\overline{\mathcal{M}}$ which directly leads to a super-linear convergence  for the new learning problem
\begin{theorem}\label{thm:transfer_learning}
Let $c_{\gamma,\rho}:=\max\{\frac{\gamma}{1-\gamma},\frac{\gamma \rho}{1-\gamma \rho^2}\}$ and $\delta' := \min \Big\{ \frac{1}{c'_1 + c'_2} \sigma_{\min}(S_{\overline{K}^*,\overline{\Sigma}^*}),\\ \frac{\rho - \|\overline{A} -\overline{B}\overline{K}^* \|}{\| \overline{B}\|} \Big\}$. If the following condition is satisfied:
\begin{eqnarray}\label{eq:close}
  \qquad \Big(\|A-\overline A \|_q +\|B-\overline B \|_q\Big) \leq \frac{\left( \frac{1}{\mu} - \frac{1}{\|S^*\|}\right)^{-1}\sigma_{\min}\left( R + \gamma \overline{B} \overline{P}_{\overline{K}^*} \overline{B}\right) (\delta')^2}{4 c_{\gamma,\rho} \Big(\left\|\mathbb{E}_{x_0\sim\mathcal{D}}[x_0 x_0^\top]+W\right\|_q + \frac{\gamma}{1-\gamma}+1\Big) \frac{\|Q\| + \|R\|}{1-\gamma \rho^2}},
  \end{eqnarray}
then $(K^*,\Sigma^*)$ is within the super-linear convergence region of environment $\overline{\mathcal{M}}$, \textit{i.e.,} 
\begin{eqnarray*}
   &&     \overline{C}(\overline{K}^{(t+1)},\overline{\Sigma}^{(t+1)}) - \overline{C}(\overline{K}^*, \overline{\Sigma}^{*})
        \leq \frac{(c'_1 +  c'_2) \cdot \left(\overline{C}(\overline{K}^{(t)}, \overline{\Sigma}^{(t)}) - \overline{C}(\overline{K}^*, \overline{\Sigma}^{*}) \right)^{1.5}}{\sigma_{\min}(S_{\overline{K}^*, \overline{\Sigma}^{*}}) \sqrt{\mu \sigma_{\min}(R+\gamma \overline{B}^\top \overline{P}_{\overline{K}^*}\overline{B})}},
    \end{eqnarray*}  for all $t \geq 0$,  
    if the initial policy follows $(\overline{K}^{(0)},\overline{\Sigma}^{(0)}) = (K^*,\Sigma^*)$ and the policy updates according to \eqref{alg: PI}.
\end{theorem}
\begin{proof} (of Theorem \ref{thm:transfer_learning}).
It is easy to verify that
\begin{equation}
    \begin{aligned}
        \left|\overline{C}(\overline{K}^*,\overline{\Sigma}^*) - \overline{C}(K^*,\Sigma^*)\right|\leq &\left|\overline{C}(\overline{K}^*,\overline{\Sigma}^*) - C(\overline{K}^*,\overline{\Sigma}^*)\right|+ \left|\overline{C}(K^*,\Sigma^*) - C(K^*,\Sigma^*)\right|.\label{eq:C_wts}
    \end{aligned}
\end{equation}
For any given policy $(K,\Sigma)$ that is admissible to both $\mathcal{M}$ and $\overline{\mathcal{M}}$, we have $J_{K, \Sigma}(x) = x^\top P_K x + q_{K,\Sigma}$ and $\overline{J}_{K, \Sigma}(x) = x^\top \overline{P}_{{K}} x + \overline{q}_{K,\Sigma}$  with some symmetric positive definite matrices $P_K,\overline{P}_{K} \in \mathbb{R}^{n\times n}$  satisfying \eqref{eqn:P_K} and \eqref{eqn:P_K2} respectively.
In addition, the constants $q_{K,\Sigma}$ takes the form of \eqref{eqn:P_K} and $\overline{q}_{K,\Sigma}\in \mathbb{R}$ takes the form of
    %      q_{K,\Sigma} =&\frac{1}{1-\gamma} \Big[ 
    %     -\frac{\tau}{2}\left(k+\log \left((2 \pi)^k \operatorname{det} \Sigma\right)\right)+ \Tr\left(\Sigma R + \gamma(\Sigma B^{\top} P_K B+ W P_K)\right) 
    % \Big],\\
      $\overline q_{K,\Sigma} =\frac{1}{1-\gamma} \Big( -\frac{\tau}{2}\left(k+\log \left((2 \pi)^k \operatorname{det} \Sigma\right)\right)+ \Tr\left(\Sigma R + \gamma(\Sigma (\overline B)^{\top} \overline P_K \overline B + W \overline P_K) \right)
    \Big).$
Note that
\begin{align*}
& P_K - \overline P_K =  \gamma(A-BK)^{\top} P_K (A-BK) - \gamma(\overline A-\overline B K)^{\top} \overline P_K (\overline A - \overline B K)\\
 & =  \gamma(A-\overline A-(B-\overline B)K)^{\top} P_K (A-BK) +  \gamma(\overline A -\overline B K)^{\top} P_K (A-\overline A -(B-\overline B)K)\\
 & \quad + \gamma(\overline A -\overline B K)^{\top} (P_K - \overline P_K) (\overline A-\overline B K). 
\end{align*}
Hence we have
%\xinyu{Is it $\|K\|$ instead of $K$?}
%\begin{eqnarray*}
$\|P_K - \overline P_K \| \leq 2 \gamma (\|A-\overline A \| + \|K\|\|B-\overline B \|) \|P_K\| \rho +\gamma \rho^2 \|P_K - \overline P_K \|,$
%\end{eqnarray*}
and therefore, since $\gamma \rho^2<1$,
\begin{eqnarray}\label{eq:P_diff}
\|P_K - \overline P_K \| \leq \frac{2\gamma\rho}{1-\gamma \rho^2} (\|A-\overline A \| + \|K\|\,\|B-\overline B \|) \|P_K\|.
\end{eqnarray}
Similarly,\vspace{-2mm}
\begin{eqnarray*}
  \overline q_{K,\Sigma} - q_{K,\Sigma} &=& \frac{\gamma}{1-\gamma} \Big(\text{Tr}\Big( (B-\overline B)^\top P_K\,B\Big) + \text{Tr}\Big( \overline{B}^\top P_K\,(B-\overline B ) \Big) \\
  &&\quad \quad \quad \quad +\text{Tr}\Big( \overline{B}^\top (P_K-\overline P_K)\,\overline B \Big)+\text{Tr}(W(P_K-\overline P_K))\Big).
\end{eqnarray*}
Recall that $x_0 \sim \mathcal{D}$ and denote $D_0 = \mathbb{E}_{x_0 \sim \mathcal{D}}[x_0 x_0^T] $. Therefore,
\begin{eqnarray}
&&\left|\overline{C}(\overline{K}^*,\overline{\Sigma}^*) - C(\overline{K}^*,\overline{\Sigma}^*)\right| \nonumber\\
&\leq& \Big|\Tr \Big((P_{\overline{K}^*} - \overline P_{\overline{K}^*})(D_0+W)\Big)\Big| + \frac{\gamma}{1-\gamma}
\left|\text{Tr}\Big( (B-\overline B)^\top P_{K^*}\,B\Big)\right|\nonumber\\
&&+\frac{\gamma}{1-\gamma}\left| \text{Tr}\Big( \overline{B}^\top P_{\overline{K}^*}\,(B-\overline B ) \Big) \right| + \frac{\gamma}{1-\gamma}\left| \text{Tr}\Big( \overline{B}^\top (P_{\overline{K}^*}-\overline P_{\overline{K}^*})\,\overline B \Big)\right|\nonumber\\
&\leq &\|P_{\overline{K}^*} - \overline P_{\overline{K}^*}\|_p \|D_0+W\|_q+ \frac{\gamma}{1-\gamma}\|B-\overline B\|_p\|P_{\overline{K}^*}\|_q (\|B\|_q+\|\overline{B}\|_q) \nonumber\\
&&+ \frac{\gamma}{1-\gamma}\| P_{\overline{K}^*}-\overline P_{\overline{K}^*}\|_p\|\overline{B}\,\overline{B}^\top\|_q, \label{eq:C_bound1}
\end{eqnarray}
where the last inequality follows from $|\text{Tr}(A B)| \leq \|A'\|_p  \, \|B\|_q$ when $1/p + 1/q = 1$.
This is a consequence of combining von Neumann's trace inequality with H\"older's inequality for Euclidean space. %yields a version of Hölder's inequality for Schatten norms for 

\eqref{eqn:P_K} suggests that, for any admissible $K$ such that $\|K\|\leq 1$, we have
$\|P_K\| \leq \gamma \rho^2 \|P_K\| +\|Q\| + \|R\|.$
%\end{eqnarray}
Hence, $\|P_K\|\leq \frac{\|Q\| + \|R\|}{1-\gamma \rho^2}$. Combining this bound with \eqref{eq:P_diff} and the fact that $\|B\|_q,\|\overline{B}\|_q,\|K^*\| \leq 1$, \eqref{eq:C_bound1} can be bounded such that 
\begin{align}
     \eqref{eq:C_bound1} &\leq \Big(\|D_0+W\|_q + \frac{\gamma}{1-\gamma}\Big)\frac{2\gamma\rho}{1-\gamma \rho^2}  \frac{\|Q\| + \|R\|}{1-\gamma \rho^2}\Big(\|A-\overline A \| +\|B-\overline B \|\Big)\notag\\
&+ \frac{2\gamma}{1-\gamma}   \frac{\|Q\| + \|R\|}{1-\gamma \rho^2}\|B-\overline{B}\|_q \notag\\
&\leq 2 c_{\gamma,\rho} \Big(\|D_0+W\|_q + \frac{\gamma}{1-\gamma}+1\Big) \frac{\|Q\| + \|R\|}{1-\gamma \rho^2}\Big(\|A-\overline A \|_q +\|B-\overline B \|_q\Big).\label{eq:C_bound3}
\end{align}
Similarly, we have
\begin{equation}\label{eq:C_bound4}
    \begin{aligned}
        &\left|\overline{C}(K^*,\Sigma^*) - C(K^*,\Sigma^*)\right|\\
 &\leq 2 c_{\gamma,\rho} \Big(\|D_0+W\|_q + \frac{\gamma}{1-\gamma}+1\Big) \frac{\|Q\| + \|R\|}{1-\gamma \rho^2}\Big(\|A-\overline A \|_q +\|B-\overline B \|_q\Big). 
    \end{aligned}
\end{equation}
Finally, plugging \eqref{eq:C_bound3}-\eqref{eq:C_bound4} into \eqref{eq:C_wts}, we have 
\begin{equation}
    \begin{aligned}
 &\left|\overline{C}(\overline{K}^*,\overline{\Sigma}^*) - \overline{C}(K^*,\Sigma^*)\right|\\
 &\leq 4 c_{\gamma,\rho} \Big(\|D_0+W\|_q + \frac{\gamma}{1-\gamma}+1\Big) \frac{\|Q\| + \|R\|}{1-\gamma \rho^2}\Big(\|A-\overline A \|_q +\|B-\overline B \|_q\Big). %\label{eq:C_bound5} 
    \end{aligned}
\end{equation}
\end{proof}

{\color{black}
\section{Model-free Extension}
\label{sec:model_free}

Model-based convergence provides a {\it foundation} for model-free analysis, as demonstrated  \cite{fazel2018global, hambly2021policy} where the  more challenging policy convergence analysis for the model-based setting is followed by a relatively routine sample-based analysis of zeroth-order gradient approximation \cite{scheinberg2000derivative,larson2019derivative,balasubramanian2018zeroth,cheng2021convergence}.
Similarly,  our analysis of local superlinear convergence and transfer learning applications within a model-based framework   can be extended to developing a model-free algorithm, for instance, Algorithm \ref{alg:modelfree}.

\begin{algorithm}\caption{Policy Gradient Estimation with Unknown Parameters}\label{alg:modelfree}
    \begin{algorithmic}
        \STATE{\textbf{Input:} $K, \Sigma$, the number of trajectories $m$, smoothing parameter $r$, dimension $D_K$ and $D_\Sigma$}
%\FOR{\(t = 0,1,2,...,T-1\)}
\STATE{Apply Cholesky decomposition to matrix $\Sigma$ to get $L$ such that
$\Sigma = L L^\top$.}
    \FOR{$i= 1, 2, \cdots, m$}
    \STATE Sample a policy $\widehat{\Sigma}_i =\widehat{L}_i (\widehat{L}_i)^\top$, where $ \vec(\widehat{L}_i) = \vec(L) + U_i$, where $U_i$ is drawn uniformly at random over
vectors in $\mathbb{R}^{D_\Sigma}$ such that $\| U_i\|_F = r.$ Simulate with policy $(K, \widehat{\Sigma}_i)$ from $x_0 \sim \mathcal{D}$ for $l$ steps. Let $\widehat{C}_i$ denote the empirical estimates:
$$
\widehat{C}_i=\sum_{t=1}^{l}\gamma^t c_t,
$$
where $c_t$ and $x_t$ are the costs and states on this trajectory.

\STATE{Sample a policy $\widehat{K}_i  =K + U_i'$, where $U_i'$ is drawn uniformly at random over matrices in $\mathbb{R}^{k \times n}$ such that $\| U_i'\|_F = r$. Simulate with policy $(\widehat{K}_i, {\Sigma}_i)$ from $x_0 \sim \mathcal{D}$ for $l$ steps and let $\widehat{C_i'}$ and $\widehat{S_i'}$ denote the empirical estimates:
$$
\widehat{C_i'} = \sum_{t=1}^{l} \gamma^t c_t', \quad \widehat{S_i'} = \sum_{t=1}^{l} \gamma^t x_t' x_t^\top.
$$}
    \ENDFOR
\RETURN{the estimates of $\nabla_{K} C(K,\Sigma), \nabla_{\Sigma}C(K,\Sigma), S_{K,\Sigma}$:
{
$$
\begin{aligned}
     \widehat{\nabla_{\vec(\Sigma)} C(K, \Sigma)} &= \left( \nabla_{\vec(L)} \vec(\widehat{\Sigma}(\widehat{L}))^\top \right)^{-1 } \cdot \left( \frac{1}{m} \sum_{i=1}^m \frac{D_\Sigma}{r^2} \widehat{C}_i U_i\right),\\ %, \quad \widehat{S} = \frac{1}{m} \sum_{i=1}^m \widehat{S}_i.
\widehat{\nabla_K C(K,\Sigma)} &= \frac{1}{m} \sum_{m=1}^l \gamma^t \frac{D_K}{r^2} \widehat{C_i'} U_i', \quad \widehat{S_{K,\Sigma}} = \frac{1}{m} \sum_{i=1}^m \widehat{S_i'}.
\end{aligned}
$$}
} % 
%\ENDFOR
\end{algorithmic}\label{alg: PolicyGradient}
\end{algorithm}

In the setting with unknown model parameters  $A, B, Q, R$, where the controller has only simulation access to the model,  
we apply a zeroth-order optimization method to approximate the gradient $\widehat{\nabla_K C(K,\Sigma)} $  and $\widehat{\nabla_{\Sigma} C(K, \Sigma)}$, as in  Algorithm \ref{alg:modelfree}.

The updating rule for  $\widehat{\nabla_K C(K,\Sigma)} $ is standard \cite{fazel2018global, hambly2021policy}. We now explain the expression for $\widehat{\nabla_{\Sigma} C(K, \Sigma)}$. By Lemma 4.5 and Theorem 4.2, $\Sigma$ is positive definite (with the time index $t$ omitted for ease of exposition). Let $L$ denote the Cholesky decomposition of $\Sigma$ such that
%\begin{align*}
   $ \Sigma = L L^\top.$
%\end{align*}
{Let $\vec (\Sigma) \in \mathbb{R}^{D_\Sigma}$ and  $\vec(L) \in \mathbb{R}^{D_\Sigma}$ denote the stacked vectors of the lower-triangular entries in
matrices $\Sigma \in \mathbb{R}^{k \times k}$ and $L \in \mathbb{R}^{k \times k}$ respectively, with $D_\Sigma :=\frac{k(k+1)}{2}$. Similarly, denote $\vec(K) \in \mathbb{R}^{D_K}$ as the stacked vectors of all entries in $K \in \mathbb{R}^{k \times n}$, with  $D_K := k \times n.$
% {\color{black}[disconnected.] Note that here, when bounding the Frobenius norm of a matrix, we usually treat the matrices ($K, \Sigma, L$ and $U$ in  Algorithm \ref{alg:modelfree}) as stacked vectors. Therefore we denote by $D_K=k \times n$ and $D_\Sigma = k \times k$ the dimensions of the corresponding vectors, transformed from the $K$ and $\Sigma$ matrices for convenience.} 
Then 
%\begin{align}\label{eq: chain_rule_sigma}
$ \nabla_{\vec(L)} C(K,\Sigma) = \nabla_{\vec(L)} \vec(\Sigma(L))^\top \cdot \nabla_{\vec(\Sigma)} C(K,\Sigma)$.
 % \\
 % \nabla_{\vec(L)} \vec(\Sigma(L))^{-1 \top}  \cdot  \nabla_{\vec(L)} C(K,\Sigma) =\nabla_{\vec(\Sigma)} C(K,\Sigma).$
%\end{align}
In Algorithm \ref{alg:modelfree}, we approximate $\nabla_{\vec(L)} C(K,\Sigma)$ with zeroth-order estimate and the above equation 
%\eqref{eq: %chain_rule_sigma} 
to get $\widehat{\nabla_{\vec(\Sigma)} C(K,\Sigma)}$. The estimate $\widehat{\nabla_{\Sigma} C(K,\Sigma)}$ can be obtained by rearranging the entries of $\widehat{\nabla_{\vec(\Sigma)} C(K,\Sigma)}$ into a matrix form.
}

}

\section{Numerical Experiments}\label{sec: numerical experiments}%\vspace{-3mm}
This section provides numerical experiments using \eqref{alg: regularized update} and \eqref{alg: PI} to illustrate the results established in Section \ref{sec:regularized_policy}, \ref{sec: policy_iteration},
and \ref{sec: transferlearning}. 

\paragraph{Setup} (1) Parameters: $A\in \mathbb{R}^{n \times n}$, $B \in \mathbb{R}^{n\times k}$, $Q \in \mathbb{S}^{n}_{+}$, $R \in \mathbb{S}^{k}_{++}$ are generated randomly. The scaling of $A$ is chosen so that $A$ is stabilizing with high probability ($\sigma_{\max}(A) < \frac{1}{\sqrt{\gamma}}$). Initialization: $K^{(0)}_{i,j} = 0.01$ for all $i, j$, $\Sigma^{(0)} = I.$ 

(2) Transfer learning setup: $\overline{A}$ and $\overline{B}$ are the state transition matrices which are generated by adding a perturbation to $A$ and $B$: $\overline{A}_{i,j} = A_{i,j} + u_{i,j}$, $\overline{B}_{i,j} = B_{i,j} + u'_{i,j}$, where $u_{i,j}$ and $u'_{i,j}$ are sampled from a uniform distribution on $[0, 10^{-3}]$. The initialization of $\overline{K}$ and $\overline{\Sigma}$ are the optimal solution $K^*$ and $\Sigma^*$ with state transition matrices $A$ and $B$.

\begin{figure}[H]
         \centering
    \begin{subfigure}[b]{0.32\textwidth}
    \includegraphics[width=4cm]{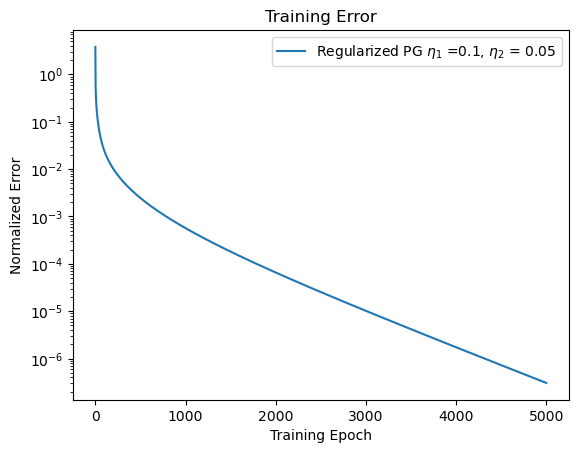}
         \caption{}
         \label{fig:reguPG}
     \end{subfigure}
     \begin{subfigure}[b]{0.32\textwidth}
    \includegraphics[width=4cm]{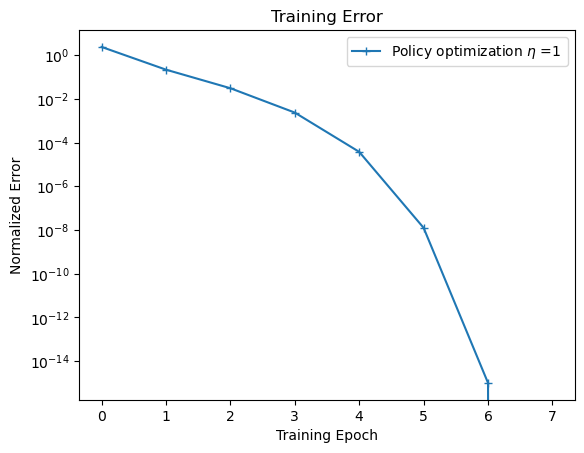}
         \caption{}
         \label{fig:unperturbed}
     \end{subfigure}
     \begin{subfigure}[b]{0.32\textwidth}
    \includegraphics[width=4cm]{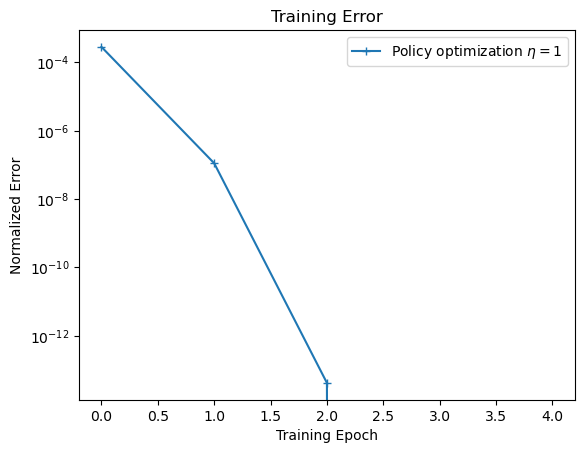}
         \caption{}
         \label{fig:perturbed}
     \end{subfigure}
    \caption{Performances of (\ref{fig:reguPG})  \eqref{alg: regularized update}; (\ref{fig:unperturbed})  \eqref{alg: PI};  (\ref{fig:perturbed}) transfer learning using \eqref{alg: PI} with $(\overline{K}^{(0)},\overline{\Sigma}^{(0)}) = (K^*,\Sigma^*)$ and state transitions $(\overline{A}, \overline{B})$. $n=400$, $k=200.$ The regularization parameter $\tau$ is chosen to be $\sigma_{\min}(R)$. }
    \label{fig:enter-label}
\end{figure}

\paragraph{Performance measure} We use normalized error to quantify the performance of a given policy ${K, \Sigma}$, \textit{i.e.,}
$
\textsc{normalized error }=\frac{C({K, \Sigma})-C\left({K}^*, \Sigma^*\right)}{C\left({K}^*,  \Sigma^*\right)},
$
where ${K}^*, \Sigma^*$ is the optimal policy defined in Theorem \ref{thm: optimal val func and policy}.
\paragraph{(Fast) Convergence} Figure \ref{fig:reguPG} shows the linear convergence of \eqref{alg: regularized update}, and Figure \ref{fig:unperturbed} shows the superior convergence rate of \eqref{alg: PI}. The normalized error falls below $10^{-14}$ within just 6 iterations, and from the third iteration, it enters a region of super-linear convergence. Figure \ref{fig:perturbed} shows the result of applying transfer learning using \eqref{alg: PI} in a perturbed environment, when the optimal policy in Figure \ref{fig:reguPG} and \ref{fig:unperturbed} serves as an initialization. Figure \ref{fig:perturbed} shows that if the process commences within a super-linear convergence region, then the error falls below $10^{-12}$ in just two epochs.

\begin{figure}[H]
    \centering
    \includegraphics[width=0.35\linewidth]{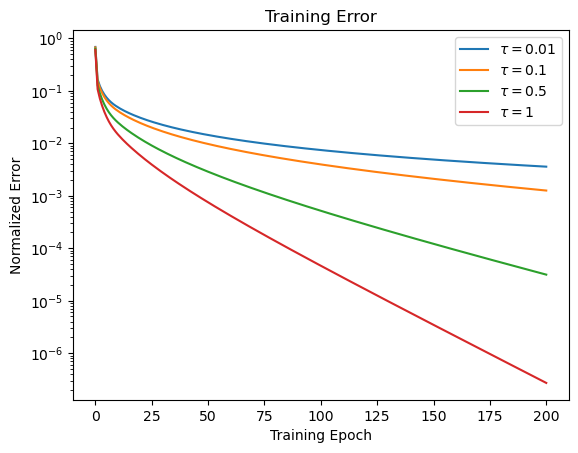}
    \includegraphics[width=0.35\linewidth]{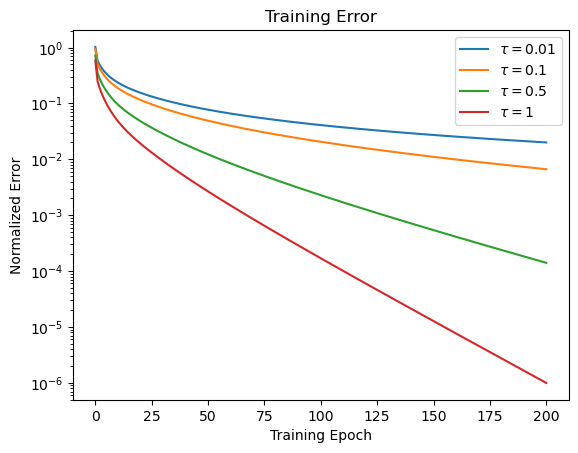}
    \caption{Regularized Policy Gradient \eqref{alg: regularized update} with different regularization parameters $\tau$. Left: $n=200, k =10;$ Right: $n=200, k =50.$}
    \label{fig:RPG_with_tau}
\end{figure}
\begin{figure}[H]
    \centering
    \includegraphics[width=0.40\linewidth]{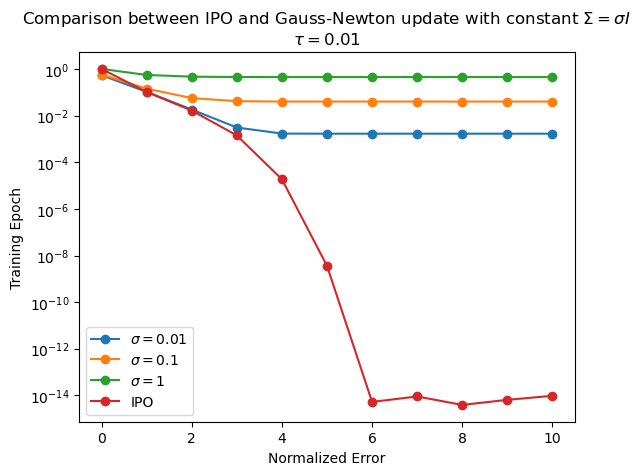}
    \includegraphics[width=0.40\linewidth]{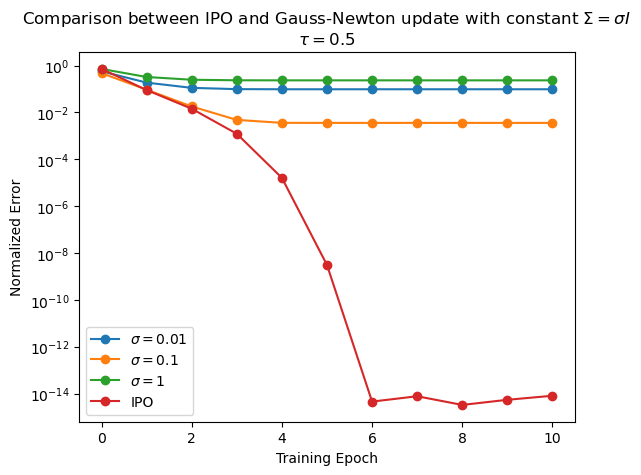}
    \caption{Comparison between Iterative Policy Optimization \eqref{alg: PI} and Gauss-Newton update on $K$ with constant covariance matrix \eqref{eq: Gauss-Newton}. $n=200, k=50$. Left: $\tau = 0.01;$ Right: $\tau = 0.5.$  %$A,B,Q,R$ are generated randomly and satisfies Assumption 1.
    }
    \label{fig:IPO_GN}
\end{figure}

{\color{black}
\paragraph{Regularization parameter $\tau$} To demonstrate that entropy regularization can accelerate convergence, we conduct experiments with \eqref{alg: regularized update} under two settings: $n=200, k=10$ and $n=200, k=50$. We run \eqref{alg: regularized update} using various values of $\tau$. Figure \ref{fig:RPG_with_tau} illustrates that, in both settings, a larger value of $\tau$ results in a faster linear convergence rate to the optimal solution of \eqref{eqn:C_K}. These results confirm that increasing the regularization parameter enhances the convergence speed, highlighting the practical benefits of entropy regularization in achieving faster optimization.

\paragraph{The importance of updating $\Sigma$}
As discussed in Section \ref{sec: policy_iteration},  updating $K$ in \eqref{alg: PI} is identical to  updating step of the Gauss-Newton algorithm in an unregularized setting. However, unlike the Gauss-Newton algorithm, \eqref{alg: PI} also updates the covariance matrix $\Sigma$ simultaneously. 
%In Figure \ref{fig:IPO_GN}, 
To see the effect of dynamic updating of  $\Sigma$, we compare the performance of \eqref{alg: PI} with  the Gauss-Newton update on $K$ with a constant covariance matrix $\Sigma$ given by:
\begin{equation}\label{eq: Gauss-Newton}
    \begin{aligned}
K^{(t+1)} &= K^{(t)} - \left(R + \gamma B^{\top} P_{K^{(t)}} B\right)^{-1} E_{K^{(t)}}, \\
\Sigma^{(t+1)} &= \sigma I,
\end{aligned}
\end{equation}
where $\sigma$ is a fixed positive scalar. Figure \ref{fig:IPO_GN} illustrates that
%the choice of covariance matrix significantly influences the convergence rate. Specifically, 
the \eqref{alg: PI} algorithm with updates of $\Sigma$ achieves a noticeably faster and superlinear convergence rate, compared to the Gauss-Newton update with a fixed $\Sigma$. This dynamic update of $\Sigma$ allows \eqref{alg: PI} to reach the optimal point of Problem \eqref{eqn:C_K} more efficiently, highlighting the importance of adapting the covariance matrix during iterations.
}
\section{Detailed Proofs} \label{sec: detailed proofs} 
\subsection{Proofs in Section \ref{sec: Problem Formulation}}\label{sec: proof of thm 1}
\subsubsection{Proof of Lemma \ref{claim: solve regularzied QP}}
% Let $\lambda$ be a Lagrangian multiplier to the constraint $\int_{\mathcal{A}} p(u)\mathrm{d}u=1.$ 
% Consider
%     \begin{align*}
%         \mathcal{L}(p, \lambda) &= \int_{\mathcal{A}} \left(u^{\top} M u + b^{\top}u + \tau \log p(u) \right) p(u) \mathrm{d}u + \lambda \left(\int_{\mathcal{A}} p(u) \mathrm{d}u - 1\right) \\
%         &=\int_{\mathcal{A}}  \left(u^{\top} M u + b^{\top}u + \tau \log p(u) + \lambda \right)p(u) \mathrm{d}u  -  \lambda\\
%         &=\int_{\mathcal{A}} {L}\left(u, p(u), \lambda\right)\mathrm{d}u - \lambda,
%     \end{align*}
% where ${L}(u,v,\lambda):= (u^{\top}Mu + b^{\top}u)v +  \tau v \log v + \lambda v.$ Since 
% $
%     \frac{\partial{L}}{\partial{v}} (u,v,\lambda) = \lambda + u^{\top}Mu + b^{\top}u +  \tau + \tau \log v,
% $
% by setting $ \frac{\partial{L}}{\partial{v}} (u,v^*,\lambda) = 0,$ we then get 
% $v^* \propto e^{-\frac{1}{\tau}(u^{\top} M u + b^{\top} u)}.$
% {\color{blue} By \cite[Chapter 8.3, Theorem 1]{luenberger1997optimization}, the strong duality holds.}
% Thus, the optimal solution $p^*$ is a multivariate Gaussian distribution with $\mathcal{N}\left(-\frac{1}{2}M^{-1}b, \frac{\tau}{2}M^{-1}\right)$.
% %\end{proof}

%{\color{purple}[I modified the proof; please double check if the new version  is correct before removing my comment]}

{Denote the domain of the decision variable as $\mathcal{X} = \{p: \mathcal{A} \mapsto [0, \infty) \},$ and the feasible set as $ \mathcal{F} = \{p: \mathcal{A} \mapsto [0, \infty) \mid \int_\mathcal{A} p(u) \mathrm{d} u = 1\}\subseteq\mathcal{X}.$ Let $f: \mathcal{X}\mapsto \mathbb{R}$ denote the objective function, i.e., $$f(p) = \mathbb{E}_{u \sim p(\cdot)}\left[u^{\top} M u+b^{\top} u+\tau \log p(u)\right].$$ 
    Let $\lambda$ be a Lagrangian multiplier to the constraint $\int_{\mathcal{A}} p(u)\mathrm{d}u=1.$ 
Consider
    \begin{align*}
        \mathcal{L}(p, \lambda) &= \int_{\mathcal{A}} \left(u^{\top} M u + b^{\top}u + \tau \log p(u) \right) p(u) \mathrm{d}u + \lambda \left(\int_{\mathcal{A}} p(u) \mathrm{d}u - 1\right) \\
        &=\int_{\mathcal{A}} {L}\left(u, p(u), \lambda\right)\mathrm{d}u - \lambda,
    \end{align*}
where ${L}(u,v,\lambda):= (u^{\top}Mu + b^{\top}u)v +  \tau v \log v + \lambda v.$ 
 Additionally,  define $g(\lambda) = \inf_{\mathcal{X}} \mathcal{L}(p, \lambda)$.
 
We now show the strong duality result: 
\begin{eqnarray}\label{eq:strong_duality}
 g(\lambda^*) = \inf_{p\in\mathcal{F}} f(p),
 \end{eqnarray}
 with $\lambda^* = \arg \max_{\lambda \in \mathbb{R}} g(\lambda).$
 
First, the weak duality result follows from
\begin{eqnarray}\label{eq:g_lambda}
    g(\lambda) = \inf_{p\in\mathcal{X}} \mathcal{L}(p,\lambda) \leq \inf_{p\in\mathcal{F}} \mathcal{L}(p,\lambda) = \inf_{p\in\mathcal{F}} f(p), \text{ for any } \lambda \in \mathbb{R}.
\end{eqnarray}
Moreover, since 
$
    \frac{\partial{L}}{\partial{v}} (u,v,\lambda) = \lambda + u^{\top}Mu + b^{\top}u +  \tau + \tau \log v,
$
for any $\lambda\in\mathbb{R}, u \in \mathcal{A}$, the minimizer $p_\lambda(u)$ of $L(u, \cdot, \lambda)$ satisfies
\begin{equation} \label{eqn: p_lambda}
    p_\lambda(u ) = \exp\left(-\frac{1}{\tau}(\lambda + u^\top M u + b^\top u) - 1\right).
\end{equation}
Therefore, by applying \eqref{eqn: p_lambda} to \eqref{eq:g_lambda}, we have
\begin{equation}
    g(\lambda) = \mathcal{L}(p_\lambda, \lambda) = - \tau \left( \exp\left(-\frac{\lambda}{\tau} - 1 \right) \cdot C + \frac{\lambda}{\tau}\right),
\end{equation}
where $C := \int_\mathcal{A} \exp\left(-\frac{1}{\tau}(u^\top M u + b^\top u)\right) \mathrm{d} u$. Direct computation yields the maximizer of $g$ in \eqref{eq:g_lambda} as  $\lambda^* =\tau \log C -\tau$. Plugging $\lambda^*$ to \eqref{eqn: p_lambda} shows $\int_\mathcal{A} p_{\lambda^*} (u) \mathrm{d} u = 1,$   implying $p_{\lambda^*} \in \mathcal{F}$ and the strong duality  \eqref{eq:strong_duality} holds. Finally, by \eqref{eq:strong_duality} and  \eqref{eqn: p_lambda}, it is clear that the optimal solution is a multivariate Gaussian distribution with $\mathcal{N}\left(-\frac{1}{2}M^{-1}b, \frac{\tau}{2}M^{-1}\right)$. }

\subsection{Proofs in Section \ref{sec: preliminary analysis}}
\subsubsection{Proof of Lemma \ref{thm: PG}}
% \begin{proof}{\bf of Lemma \ref{thm: PG}}
$\nabla_\Sigma C(K,\Sigma)$ in \eqref{eqn:grad_K} can be checked by direct gradient calculation.
 To verify $\nabla_K C(K,\Sigma)$ in \eqref{eqn:grad_K}, first %fix $K \in \mathbb{R}^{k\times n}$ and 
 define 
    $f: \mathcal{X}\times \mathbb{R}^{k\times n} \rightarrow \mathbb{R}$ 
by
    $
        f(y, K) :=  y^{\top} P_K y, \forall y \in \mathcal{X}
    $
and we aim to find the gradient of $f$ with respect to $K$.
For any $ y\in\mathcal{X},$ by the Riccati equation for $P_K$ in \eqref{eqn:P_K} we have %\vspace{-3mm}
    \begin{align*}
        f(y, K) &=  y^{\top} \left( \gamma(A-BK)^{\top} P_K (A-BK) + Q+ K^{\top}RK\right) y\\
        &= \gamma f\left((A-BK)y, K\right) + y^{\top} \left(Q+ K^{\top} R K \right)y. 
    \end{align*}
$\nabla_K f((A-BK)y, K)$ has two terms: one with respect to the input $(A-BK)y$ and one with respect to $K$ in the subscript of $P_K$. This implies \vspace{-1mm}
\begin{align}
    &\nabla_K f(y, K) = 2 \left(-\gamma B^{\top} P_K(A-BK) + R K\right) y y^{\top}+ \gamma \nabla_K f(y', K)|_{y' = (A-BK)y} \label{eqn:grad_K f} \\
     &= 2 \left(-\gamma B^{\top} P_K(A-BK) + R K\right)  \sum_{i=0}^\infty \gamma^i (A-BK)^i y y^{\top} (A^{\top}-K^{\top}B^{\top})^i. \nonumber 
\end{align}
Since $C(K,\Sigma) = \mathbb{E}_{x_0 \sim\mathcal{D}} [x_0 P_K x_0] + q_{K,\Sigma}$ with $P_K$ and $q_{K, \Sigma}$ satisfying \eqref{eqn:P_K}, 
then the gradient of $C(K, \Sigma)$ with respect to $K$ is
$$\begin{aligned}
    \nabla_K C(K, \Sigma) =& \mathbb{E}[ \nabla_K f(x_0, K)  ] + \nabla_K q_{K,\Sigma}\\
    %=&  \mathbb{E}_{x_0 \sim \mathcal{D}}\left[ 2 \left(-\gamma B^{\top} P_K(A-BK) + R K\right) x_0 x_0^{\top}+ \gamma \nabla_K f(\hat{x}_1)|_{\hat{x}_1 = (A-BK)x_0}\right] \\&+ \frac{\gamma}{1-\gamma} \left(  \nabla_K \left( \Tr(\Sigma B^\top P_K B ) \right) +  \nabla_K \left(\text{Tr}\left(W P_K\right)\right) \right)\\
    \overset{(a)}{=}& \mathbb{E}\Big[ 2 \left(-\gamma B^{\top} P_K(A-BK) + R K\right) x_0 x_0^{\top} 
     + \gamma \nabla_K f(\hat{x}_1, K)|_{\hat{x}_1 = (A-BK)x_0}\\
        & + \sum_{t=0}\gamma^{t+1}  \left( \nabla_K( \Sigma B^\top P_K B ) + \nabla_K(  w_t^{\top} P_K w_t) \right) \Big]\\
    \overset{(b)}{=}& \mathbb{E}\Big[ 2 \left(-\gamma B^{\top} P_K(A-BK) + R K\right) x_0 x_0^{\top} 
    + \gamma \nabla_K f({x}_1, K)|_{{x}_1 = A+Bu_0 + w_0} \\
        & + \sum_{t=1}\gamma^{t+1}  \left( \nabla_K( \Sigma B^\top P_K B ) + \nabla_K(  w_t^{\top} P_K w_t) \right) \Big]\\
     \overset{(c)}{=}&  2 \left(-\gamma B^{\top} P_K(A-BK) + R K\right) \mathbb{E} \left[ \sum_{i=0}^\infty \gamma^i x_i x_i^{\top}\right],
\end{aligned}$$
where $(a)$ follows from applying \eqref{eqn:grad_K f} with $y=x_0$ and taking the gradient of $q_{K,\Sigma}$ in \eqref{eqn:P_K} with respect to $K$; $(b)$ follows from \vspace{-1.5mm}
\begin{equation*}
    \mathbb{E}_{x_0, u_0, w_0}[f(A x_0 +B u_0 +w_0), K] = \mathbb{E}_{x_0, w_0} [f\left( (A-BK)x_0, K\right) + w_0^\top P_K w_0] + \Sigma B^\top P_K B.
    \vspace{-1.5mm}
\end{equation*}
Using recursion to get $(c)$. 

\subsubsection{Proof of Lemma \ref{lemma:grad_dom}}\label{sec: proof of grad_dom}
For a given policy $\pi_{K, \Sigma}(\cdot | x) = \mathcal{N}(-Kx, \Sigma)$ with parameter $K$ and $\Sigma$, we define the state-action value function (also known as $Q$-function) $Q_{K, \Sigma}: \mathcal{X} \times \mathcal{A} \rightarrow \mathbb{R}$ as the cost of the policy starting with $x_0 = x$, taking a fixed action $u_0=u$ and then proceeding with $\pi_{K,\Sigma}$. The $Q$-function is related to the value function $J_{K,\Sigma}$ defined in \eqref{eqn: obj initial} as
\begin{equation}
    \begin{aligned}
     Q_{K,\Sigma}(x,u) = & x^{\top} Q x + u^{\top} R u + \tau  \pi_{K,\Sigma}(u|x)+\gamma 
    \mathbb{E}\left[ J_{K,\Sigma}(Ax + Bu + w)
    \right] \label{eqn:Q updates},
\end{aligned}
\end{equation}
for any $(x, u)\in \mathcal{X} \times \mathcal{A}$.
By definition of the $Q$-function, we also have the relationship
    $ J_{K,\Sigma}(x) =\mathbb{E}_{u \sim\pi(\cdot|x)} \left[
       Q_{K, \Sigma}(x,u)
    \right], \, \forall x\in \mathcal{X}.$
We then introduce the advantage function $A_{K,\Sigma} :\mathcal{X} \times \mathcal{A}\rightarrow \mathbb{R}$ of the policy $\pi$:\vspace{-1mm}
\begin{equation}\label{eqn: A}
   A_{K,\Sigma}(x,u) = Q_{K,\Sigma}(x,u) - J_{K,\Sigma}(x), \vspace{-1mm}
\end{equation}
which reflects the gain one can harvest by executing control $u$ instead of following
the policy $\pi_{K, \Sigma}$ in state $x$.

With the notations of the Q-function in \eqref{eqn:Q updates} and the advantage function in \eqref{eqn: A}, we first provide a convenient form for the difference of the objective functions with respect to two different policies in Lemma \ref{lemma:cost_diff}. It takes a similar form as \cite[Lemma 10]{fazel2018global}, but with an additional term on $\Sigma$.  This will be used in the proof of Lemma \ref{lemma:grad_dom}.
\begin{lemma}[Cost difference]\label{lemma:cost_diff}Suppose policies $\pi$ and $\pi'$ are in form of \eqref{eqn: pi_simple} with parameters $(K', \Sigma')$ and $(K, \Sigma).$ Let $\{ x_t'\}_{t=0}^\infty$ and $\{ u_t'\}_{t=0}^\infty$ be state and action sequences generated by $\pi'$ with noise sequence $\{ w_t'\}_{t=0}^\infty$ (i.i.d with mean $0$ and covariance $W$), \textit{i.e.,} $x_{t+1}' =  A x_t'+B u_t' + w_t'$. Then for any $x\in \mathcal{X},$\vspace{-1mm}
\begin{equation} 
    J_{K', \Sigma'}(x) - J_{K, \Sigma}(x) = 
    \mathbb{E}_{\pi'}
    \left[\left.\sum_{t=0}^\infty \gamma^t
          A_{K, \Sigma}(x_t', u_t') \right| x_0' = x
    \right],\vspace{-1mm}
\end{equation}
where the expectation is taken over $u_t' \sim \pi'(\cdot|x_t')$ and $w_t$ for all $t = 0, 1, 2, \cdots$.

The expected advantage for any $x \in \mathcal{X}$ by taking expectation over $u \sim \pi'(\cdot |x)$ is:
\begin{equation}\label{eqn:A K Sigma}
\begin{aligned}
    &\mathbb{E}_{\pi'}\left[A_{K, \Sigma}(x, u)\right]=
 x^{\top} (K'-K)^{\top} (R + \gamma B^{\top} P_K B) (K'-K) x \\
        & \quad + 2 x^{\top} (K'-K)^{\top} \left[ (R + \gamma B^{\top} P_K B) K - \gamma B^{\top} P_KA \right] x + (1-\gamma) (f_K(\Sigma) - f_K(\Sigma')).
\end{aligned}
\end{equation}
\end{lemma}
\begin{proof} (of Lemma \ref{lemma:cost_diff}).
    Let $\{c_t'\}_{t=0}^\infty$ be the cost sequences generated by $\pi',$ \textit{i.e.,} 
    $c_t' = x_t'^{\top} Q x_t' + u_t'^{\top} R u_t' + \tau \log \pi'(u_t' | x_t'),$
    and $J_{K', \Sigma'}(x)  = \mathbb{E}_{ \pi'}\left[\sum_{t=0}^\infty \gamma^t c_t'\right]$. Then we have
    \begin{align*}
         &J_{K', \Sigma'}(x) - J_{K, \Sigma}(x) 
         = \mathbb{E}_{ \pi'}
         \left[\sum_{t=0}^\infty \gamma^t c_t'\right] -  J_{K, \Sigma}(x) \\
         =&\mathbb{E}_{ \pi'}
         \left[\sum_{t=0}^\infty 
         \gamma^t \left(c_t' - J_{K,\Sigma}(x_t') +J_{K,\Sigma}(x_t') \right)\right] -  J_{K, \Sigma}(x) \\
         \overset{(a)}{=}&\mathbb{E}_{ \pi'}
         \left[\sum_{t=0}^\infty 
         \gamma^t 
         \left(c_t'  -J_{K,\Sigma}(x_t') \right)
         + \sum_{t=1}^\infty  \gamma^t J_{K,\Sigma}(x_t')
         \right] \\
         =&\mathbb{E}_{ \pi'}
          \left[\sum_{t=0}^\infty \gamma^t 
            \left(c_t' + \gamma J_{K,\Sigma}(x_{t+1}') -J_{K,\Sigma}(x_t') 
          \right)\right] \\
       = &\mathbb{E}_{ \pi'}
          \left[\sum_{t=0}^\infty \gamma^t 
            \left(x_t'^{\top} Q x_t' + u_t'^{\top} R u_t' + \tau \log \pi'(u_t' | x_t')
            + \gamma J_{K,\Sigma}(x_{t+1}') -J_{K,\Sigma}(x_t') 
          \right)\right],
          \end{align*}
          where $(a)$ follows from $x_0'=x$. %and the expectations are taken over ${u_t'\sim \pi'(\cdot | x_t')}$ for all $t$ and the noises $\{ w_t'\}_{t=0}^\infty$.
        By \eqref{eqn:Q updates} and \eqref{eqn: A}, we can continue the above calculations by
          \begin{align*}
           J_{K', \Sigma'}(x) - J_{K, \Sigma}(x) &=\mathbb{E}_{ \pi'}
          \left[\sum_{t=0}^\infty \gamma^t
          \left(
          Q_{K, \Sigma}(x_t', u_t') - J_{K,\Sigma}(x_t')
          \right)
          \right]=\mathbb{E}_{ \pi'}
          \left[\sum_{t=0}^\infty \gamma^t
          A_{K, \Sigma}(x_t', u_t')
          \right].
    \end{align*}
To derive the expectation of the advantage function, for any $x$ in $\mathcal{X}$, take expectation over $u \sim \pi'(\cdot |x):=\mathcal{N}(-K'x, \Sigma')$  to get:
\begin{align*}
\mathbb{E}_{ \pi'} &
 \left[  A_{K, \Sigma}(x, u)\right] 
    = \mathbb{E}_{ \pi'}[ Q_{K,\Sigma}(x,u) ] 
        - J_{K,\Sigma}(x) \\
    %= &x^{\top} (Q + K'^{\top} R K' )x + \text{Tr}(\Sigma' R) -\frac{\tau}{2}\left(k+\log \left((2 \pi)^k \operatorname{det} \Sigma'\right)\right)\\
    %    &+ \gamma \mathbb{E} \left[J_{K,\Sigma}(Ax + Bu + w) \right]  - J_{K,\Sigma}(x) \\
    =  &x^{\top} (Q + K'^{\top} R K' )x + \text{Tr}(\Sigma' R) -\frac{\tau}{2}\left(k+\log \left((2 \pi)^k \operatorname{det} \Sigma'\right)\right)\\\
        &+ \gamma  \mathbb{E}
        \left[ (Ax + Bu +w)^{\top} P_{K} (Ax + Bu +w) + q_{K,\Sigma} \right] - x^{\top} P_K x - q_{K,\Sigma}\\
    = &x^{\top} (Q + K'^{\top} R K' ) x +  \text{Tr}(\Sigma' R)  -\frac{\tau}{2}\left(k+\log \left((2 \pi)^k \operatorname{det} \Sigma'\right)\right) - x^{\top} P_K x - (1-\gamma) q_{K,\Sigma}\\
        &+ \gamma \left[ 
        x^{\top} (A - BK')^{\top} P_K (A-BK') x   + \text{Tr}(\Sigma' B^{\top} P_K B) + \text{Tr}(W P_K)
    \right]. 
\end{align*}
Plug \eqref{eqn:P_K} into the above equation, then use \eqref{eqn: f} and write $K' = K + K' - K$ to get 
\begin{align*}
      \mathbb{E}_{\pi'}
&\left[A_{K, \Sigma}(x, u)\right] \\
    %= & x^{\top} (Q + K'^{\top} R K') x +  \text{Tr}(\Sigma' R)  -\frac{\tau}{2}\left(k+\log \left((2 \pi)^k \operatorname{det} \Sigma'\right)\right)\\
    %    &+ \gamma \left[  x^{\top} (A - BK')^{\top} P_K (A-BK') x   + \text{Tr}(\Sigma' B^{\top} P_K B) + \text{Tr}(W P_K) \right] -x^{\top} P_K x \\
    % & - \text{Tr}(\Sigma R) +\frac{\tau}{2}\left(k+\log \left((2 \pi)^k \operatorname{det} \Sigma\right)\right) -\gamma \text{Tr}(\Sigma B^{\top} P_K B) - \gamma\text{Tr}(W P_K) \\
     =&  (1-\gamma) (f_K(\Sigma) - f_K(\Sigma')) + x^{\top} \left(Q + (K + K' - K)^{\top} R  (K + K' - K )\right)x \\
        &+ \gamma x^{\top} (A - BK - B(K'-K))^{\top} P_K (A - BK - B(K'-K))x - x^{\top} P_K x \\
        % &+ \text{Tr}\left((\Sigma' - \Sigma) (R + \gamma B^{\top} P_K B) \right) +\frac{\tau}{2}\left(\log \operatorname{det} \Sigma  - \log\det \Sigma'\right) \\
    =&   (1-\gamma) (f_K(\Sigma) - f_K(\Sigma')) + x^{\top} \left(Q + K^{\top} R K 
        + \gamma(A-BK)^{\top} P_K (A-BK) \right)x  \\
        & + 2 x^{\top} (K'-K)^{\top} \left[ (R + \gamma B^{\top} P_K B) K - \gamma B^{\top} P_KA \right] x \\
        &+ x^{\top} (K'-K)^{\top} (R + \gamma B^{\top} P_K B) (K'-K) x - x^{\top} P_K x.
        % &+ \Tr\left((\Sigma' - \Sigma) (R + \gamma B^{\top} P_K B) \right)+\frac{\tau}{2}\left(\log \operatorname{det} \Sigma  - \log\det \Sigma'\right)
\end{align*}       
Plugging \eqref{eqn:P_K} finishes the proof to verify \eqref{eqn:A K Sigma}.
\end{proof}

\begin{proof} (of Lemma \ref{lemma:grad_dom}).
By Lemma \ref{lemma:cost_diff}, for any $\pi$ and $\pi'$ in form of \eqref{eqn: pi_simple} with parameter $(K,\Sigma)$ and $(K', \Sigma')$ respectively, and for any $x \in \mathcal{X}$,\vspace{-1mm}
    \begin{align}
    \mathbb{E}_{ \pi'} [ A_{K, \Sigma}&(x, u)]  =   (1-\gamma) (f_K(\Sigma) - f_K(\Sigma')) -\Tr\left(x x^{\top} E_K^{\top}\left(R+\gamma B^{\top} P_K B\right)^{-1} E_K\right) \notag  \\ 
        & \quad + \Tr\Bigl[
            xx^{\top}\left(K^{\prime}-K+\left(R+\gamma B^{\top} P_K       B\right)^{-1} E_K\right)^{\top}\left(R+\gamma B^{\top} P_K      B\right) \notag   \\ 
            &\quad \quad  \cdot\left(K^{\prime}-K+(R+\gamma B^{\top} P_K B)^{-1} E_K\right)\Bigr]  \notag \\ 
        &\geq  -\Tr\left(x x^{\top}     
            E_K^{\top}\left(R+\gamma B^{\top} P_K B\right)^{-1} E_K\right) + (1-\gamma) \left(f_K\left(\Sigma\right) - f_K\left(\Sigma'\right)\right), \label{eqn:A_K, Sigma tmp}
    \end{align}
 with equality when $K - (R+\gamma B^{\top} P_K B)^{-1} E_K = K'$ holds.
Let $\{x_t^*\}_{t=0}^\infty$ and $\{u_t^*\}_{t=0}^\infty$ be the state and control sequences generated under $\pi^*(\cdot|x) = \mathcal{N}(-K^*x, \Sigma^*)$ with noise sequence $\{w_t^*\}_{t=0}^\infty$ with mean $0$ and covariance $W$. Apply Lemma \ref{lemma:cost_diff} and \eqref{eqn:A_K, Sigma tmp} with $\pi$ and $\pi^*$ to get:\vspace{-2mm}
\begin{align}%\label{eqn: C minus C*}
    &C(K, \Sigma) - C(K^*, \Sigma^*) 
    = - \mathbb{E}_{\pi^*}
          \Bigg[\sum_{t=0}^\infty \gamma^t
          A_{K, \Sigma}(x_t^*, u_t^*) 
          \Bigg] \notag \\
    &\leq \Tr\left( 
        S_{K^*,\Sigma^*} E_K^{\top}\left(R+\gamma B^{\top} P_K B\right)^{-1} E_K\right)  + f_K(\Sigma^*) - f_K(\Sigma), \label{eqn: K part in C minus C *} 
\end{align}
where $f_K$ is defined in \eqref{eqn: f}. To analyze the first term in \eqref{eqn: K part in C minus C *}, note that
\begin{equation}\label{eqn: grad K part}
    \begin{aligned}
    &\Tr\left( 
        S_{K^*,\Sigma^*} E_K^{\top}\left(R+\gamma B^{\top} P_K B\right)^{-1} E_K\right)
    \leq  \frac{\| S_{K^*,\Sigma^*}\|}{\sigma_{\min}(R)}\operatorname{Tr}(E_K^\top E_K)\\
    &\leq  \frac{\| S_{K^*,\Sigma^*} \|}{\mu^2\sigma_{\min}(R)} \Tr\left(
        S_{K,\Sigma} E_K^\top E_K S_{K,\Sigma} \right) = \frac{\| S_{K^*,\Sigma^*} \|}{4\mu^2\sigma_{\min}(R)} \Tr(\nabla_K^\top C(K, \Sigma) \nabla_K C(K, \Sigma)),\end{aligned}
\end{equation}
where the last equation follows from \eqref{eqn:grad_K}.
% To bound $\operatorname{Tr}(E_K^\top E_K)$ with $\|\nabla_K C(K, \Sigma) \|_F,$ apply \eqref{eqn:grad_K} to get
% $$\operatorname{Tr}(\nabla_K^\top C(K, \Sigma) \nabla_K C(K, \Sigma)) = 4 \operatorname{Tr}\left(
%         S_{K,\Sigma} E_K^\top E_K S_{K,\Sigma} 
%         \right)
%     \geq 4\mu^2 \operatorname{Tr}(E_K^\top E_K).$$
% Plug the above result into \eqref{eqn: grad K part} and use \eqref{eqn: f Sigma * - f Sigma} to get \eqref{eqn: grad dom bound}.
% \begin{equation}
% \begin{aligned}
%     C(K, \Sigma) - C(K^*, \Sigma^*) \leq \frac{\| S_{K^*,\Sigma^*}\| }{4\mu^2\sigma_{\min}(R)}\|\nabla_K C(K, \Sigma)\|_F^2 +  \frac{1-\gamma}{\sigma_{\min}(R)} \| \nabla_\Sigma C(K,\Sigma)\|_F^2. 
% \end{aligned}
% \end{equation}

To analyze the second terms in \eqref{eqn: K part in C minus C *},
note that $f_K$ is a concave function, thus we can find its maximizer $\Sigma_K^*$ by taking the gradient of $f_K$ and setting it to $0$, \textit{i.e.,} $\Sigma_K^* = \frac{\tau}{2}(R+ \gamma B^\top P_K B)^{-1}.$
Thus,
\begin{equation}\label{eqn: f Sigma * - f Sigma}
    \begin{aligned}
   &  f_K(\Sigma^*) - f_K(\Sigma) \leq f_K(\Sigma_K^*) - f_K(\Sigma)  \overset{(a)}{\leq} \operatorname{Tr}\left( \nabla f_K(\Sigma)^\top (\Sigma_K^* - \Sigma) \right)  \\
     % & = \operatorname{Tr}\left(- \nabla_\Sigma C(K, \Sigma) \cdot \left(\frac{\tau}{2} (R+\gamma B^\top P_K B)^{-1} - \Sigma\right)\right) \\
     & = \operatorname{Tr}\left( \nabla_\Sigma C(K, \Sigma) \cdot ( R+\gamma B^\top P_K B)^{-1} \left( ( R+\gamma B^\top P_K B)- \frac{\tau}{2} \Sigma^{-1}\right) \Sigma  \right)  \\
    & \overset{(b)}{\leq}  (1-\gamma)  \cdot \| ( R+\gamma B^\top P_K B)^{-1}\| \cdot \| \nabla_\Sigma C(K, \Sigma)\|_F^2 \leq \frac{(1-\gamma)  \| \nabla_\Sigma C(K, \Sigma)\|_F^2}{\sigma_{\min}(R)}  ,
    \end{aligned}
\end{equation}
where $(a)$ follows from the first order concavity condition for $f_K$ and $(b)$ is from $0 \prec \Sigma \preceq I.$ Plug \eqref{eqn: grad K part} and \eqref{eqn: f Sigma * - f Sigma} into \eqref{eqn: K part in C minus C *} to get \eqref{eqn: grad dom bound}.

For the lower bound, consider $K' = K - (R+\gamma B^{\top} P_K B)^{-1} E_K$ and $\Sigma' = \Sigma$ where equality holds in \eqref{eqn:A_K, Sigma tmp}. Let $\{x_t'\}_{t=0}^\infty$, $\{u_t'\}_{t=0}^\infty$ be the sequence generated with $K', \Sigma'$. By $C(K^*, \Sigma^*) \leq C(K',\Sigma'),$ we have
\begin{equation}
    \begin{aligned}
    & C(K,\Sigma) - C(K^*, \Sigma^*) 
    \geq C(K, \Sigma) - C(K', \Sigma') = - \mathbb{E}\left[ \sum_{t=0}^\infty \gamma^t A_{K, \Sigma} (x_t', u_t')\right] \\
    &= \Tr\left( S_{K',\Sigma'} E_K^{\top}\left(R+\gamma B^{\top} P_K B\right)^{-1} E_K\right)
    % &=  \sum_{t=0}^\infty  \gamma^t \mathbb{E}\left[ \operatorname{Tr}\left(x_t' (x_t')^{\top}
    %         E_K^{\top}\left(R+\gamma B^{\top} P_K B\right)^{-1} E_K\right) \right]\\
    % \geq \frac{\sigma_{\min} 
    %     (S_{K',\Sigma'})}{\| R+\gamma B^{\top} P_K     B\|}
    %     \operatorname{Tr}(E_K^\top E_K)
        \geq \frac{\mu}{\| R+\gamma B^{\top} P_K     B\|}
        \operatorname{Tr}(E_K^\top E_K).
\end{aligned}%\vspace{-1mm}
\end{equation}
\end{proof}
\subsubsection{Proof of Lemma \ref{lemma:almost_smooth}}\label{sec: proof of almost smooth}
Lemma \ref{lemma: smooth of fK} shows that  the cost objective is smooth in $\Sigma$ when utilizing entropy regularization,  given that $\Sigma$ is bounded.
\begin{lemma}[Smoothness of $f_K$ \eqref{eqn: f}]\label{lemma: smooth of fK}Let $K \in \mathbb{R}^{k \times n}$ be given and let $f_K$ be defined in \eqref{eqn: f}. Fix $0 < a \leq 1$.
For any symmetric positive definite matrices $X \in \mathbb{R}^{k\times k}$ and $Y \in \mathbb{R}^{k\times k}$ satisfying $aI \preceq X \preceq I$ and $aI \preceq Y \preceq I$, %the following inequality holds
    \begin{align*}
        \left|  f_K(X) - f_K(Y) + \operatorname{Tr}\left( \nabla f_K(X)^\top (Y-X)\right) \right| \leq M_{a}\Tr\left( (X^{-1} Y - I)^2\right),
    \end{align*}
    where $M_{a} \in \mathbb{R}$ is defined in Lemma \ref{lemma:almost_smooth}, and $M_{a} \geq \frac{\tau}{4(1-\gamma)}$.
\end{lemma}

\begin{proof} (of Lemma
\ref{lemma: smooth of fK}).
    Fix symmetric positive definite matrices $X$ and $Y$ satisfying $aI \preceq X \preceq I$ and $aI \preceq Y \preceq I$. Then  $f_K$ being concave implies 
    $
    f_K(X) - f_K(Y) + \operatorname{Tr}\left( \nabla f_K(X)^\top (Y-X)\right)  \geq 0.
    $
    To find an upper bound, observe that
    \begin{equation}\label{eqn: smooth of fk part 1}
        \begin{aligned}
          & f_K(X) - f_K(Y) + \operatorname{Tr}\left( \nabla f_K(X)^\top (Y-X)\right)   \\
          =& \frac{\tau}{2(1-\gamma)} \left( \log \operatorname{det}(X) - \log \operatorname{det}(Y) 
          + \operatorname{Tr}\left( 
           X^{-1}(Y-X)  \right) \right).
    \end{aligned}
    \end{equation}
  Since $X \succ 0$ and $Y \succ 0,$ then all eigenvalues of $X^{-1}Y$ are real and positive and $\sigma_{\min}(X^{-1} Y ) \geq  \sigma_{\min}(X^{-1}) \sigma_{\min}(Y) \geq a.$ 
  
  Now let us show that there exists $m \in \mathbb{R}^+$ (independent of $K$) such that for any $Z \in \mathbb{R}^{k \times k}$ with real positive eigenvalues $a\leq z_1 \leq \cdots \leq z_k$, the following holds:
  \begin{equation}\label{eqn: log Z}
        -\log(\det(Z)) +  \Tr(Z-I) \leq m \Tr((Z-I)^2).
  \end{equation}
  Note that \eqref{eqn: log Z} is equivalent to 
      $\sum_{i=1}^k - \log(z_i) + z_i - 1 \leq m \sum_{i=1}^k
 (z_i - 1)^2.$
With elementary algebra, one can verify that when $m := {
(-\log\left(a\right) +  a-1)}\cdot{\left(a-1\right)^{-2}},$ it holds that
$
-\log(z) + z- 1 \leq m(z-1)^2
$
for all $z \geq a$ and such an $m$ satisfies $m \geq \frac{1}{2}$.
Therefore, \eqref{eqn: log Z} holds. Combining \eqref{eqn: smooth of fk part 1}
and \eqref{eqn: log Z} with $Z =  X^{-1}Y$, we see 
$f_K(X) - f_K(Y) +\Tr\left( \nabla f_K(X)^\top (Y-X)\right)\leq  \frac{\tau m}{2(1-\gamma)} \Tr\left( ( X^{-1}Y - I)^2\right).$ %\\ &\leq \frac{\tau m}{2(1-\gamma)} \sigma_{\min}(X^{-1})^2 \|X - Y \|_F^2 \leq \frac{\tau m}{2a^2(1-\gamma)} \| X- Y\|_F^2.$
\end{proof}
\begin{proof} (of Lemma \ref{lemma:almost_smooth}).
The first equality immediately results from \eqref{eqn:A K Sigma} in Lemma \ref{lemma:cost_diff}.
%\rx{what do you mean by taking an expectation?}. 
The last inequality follows directly from Lemma \ref{lemma: smooth of fK}.
\end{proof}
\subsection{Proofs in Section \ref{sec:regularized_policy}}
\subsubsection{Proofs of Lemma \ref{lemma: bound of update Sigma}} \label{sec: proof of lemma: bound of update Sigma} 
To ease the exposition, let $\eta$ denote $\eta_2$.
     The proof is composed of two steps. First, one can show \vspace{-1mm}\begin{equation}\label{eqn: Sigma bound in PG}
        aI \preceq  \Sigma -  \eta(1-\gamma)^{-1}( R - \frac{\tau}{2}\Sigma^{-1} + \gamma B^\top P_K B) \preceq I.\vspace{-1mm}
    \end{equation}
     Let $g: \mathbb{R}^+ \rightarrow \mathbb{R}$ be a function such that $g(x) = x + \frac{{\eta}\tau}{2(1-\gamma)x}.$ Thus, $g$ monotonically increases on $\left[\sqrt{\frac{{\eta} \tau }{2(1-\gamma)}}, \infty\right)$.
    Since 
    $ \sqrt{\frac{{\eta} \tau }{2(1-\gamma)}}  \leq a \leq \frac{\sigma_{\min}(R)}{\| R + \gamma B^\top P_K B\|} \leq 1,$ then \vspace{-0.5mm}
    \begin{equation*}
       \begin{aligned}
        &\Sigma +  \frac{{\eta}\tau}{2(1-\gamma)} \Sigma^{-1} - \frac{\eta(R+\gamma B^\top P_K B)}{1-\gamma}  \succeq \left(a +  \frac{{\eta}\tau}{2(1-\gamma)a} \right)I - \frac{\eta (R+\gamma B^\top P_K B)}{1-\gamma}  \\
        &\succeq \left( a + \frac{\eta \| {R} + \gamma B^\top P_K B\|}{1-\gamma} \right) I - \frac{\eta(R+\gamma B^\top P_K B)}{1-\gamma}  
        \succeq aI, \text{ and }\\
        &\Sigma +  \frac{{\eta}\tau}{2(1-\gamma)} \Sigma^{-1} -\frac{\eta  (R+\gamma B^\top P_K B)}{1-\gamma}  \preceq \left(1 +  \frac{{\eta}\tau}{2(1-\gamma)} \right)I - \frac{\eta (R+\gamma B^\top P_K B)}{1-\gamma} \\
        & \preceq \left( 1 + \frac{\eta \sigma_{\min}(R)}{1-\gamma}  \right) I - \frac{\eta(R+\gamma B^\top P_K B)}{1-\gamma} 
        \preceq I, \text{ hence \eqref{eqn: Sigma bound in PG}.}
    \end{aligned} \vspace{-0.5mm}
    \end{equation*}
    Second, one can show $aI \preceq \Sigma' \preceq I:$ observe that \eqref{eqn: Sigma bound in PG} is equivalent to 
     $
    aI  - \Sigma \preceq -\frac{\eta}{1-\gamma} (R+\gamma B^\top P_K B -\frac{\tau}{2}\Sigma^{-1}) \preceq  I  - \Sigma.
    $
    Then by multiplying $\Sigma$ to both sides then adding a $\Sigma$ to each term, we have \vspace{-1mm}
    \begin{equation*} 
        a\Sigma^2 - \Sigma^3 +\Sigma  \preceq \Sigma-\frac{\eta}{1-\gamma}  \Sigma (R+\gamma B^\top P_K B-\frac{\tau}{2}\Sigma^{-1}) \Sigma  \preceq \Sigma^2- \Sigma^3 +\Sigma. \vspace{-1mm}
    \end{equation*}
    With  $aI-\Sigma \preceq 0$,  $I-\Sigma \succeq 0$ and $\Sigma \preceq I$, it holds that
    %\begin{align*}
      $ aI - \Sigma \preceq (a I - \Sigma) \Sigma^2,$ and
        $I - \Sigma \succeq ( I - \Sigma) \Sigma^2.$
    %\end{align*}
   This implies
    $
  aI \preceq a\Sigma^2 - \Sigma^3 +\Sigma \preceq \Sigma' \preceq \Sigma^2 - \Sigma^3 +\Sigma \preceq I.
    $

\subsection{Proof of Lemma \ref{lemma:PG contraction}}

% \begin{proof}{\bf of Lemma \ref{lemma:PG contraction}}
     For ease of  exposition, write $S = S_{K,\Sigma}, $  $S' = S_{K',\Sigma'} $, and $S^* = S_{K^*,\Sigma^*} $. Let $M_{a}$ be defined in the same way as in Lemma \ref{lemma:almost_smooth}. Let $f_K$ be defined as \eqref{eqn: f}. Then Lemma \ref{lemma:almost_smooth} implies 
     \begin{equation}\vspace{-1mm}
     \begin{aligned}
          C(K', \Sigma') - C(K, \Sigma)&= \operatorname{Tr}\left(S' (K'-K)^\top (R +\gamma B^\top P_K B)(K'-K)\right)\\
          &\quad+ 2 \operatorname{Tr}\left(S' (K'-K)^\top E_K \right)  + f_K(\Sigma) - f_K(\Sigma'). \label{eqn: Cprime minus C}
     \end{aligned}\vspace{-1mm}
    \end{equation}
By \eqref{alg: regularized update},\vspace{-0.5mm} %Lemma \ref{lemma:grad_dom},
\begin{equation}\label{eqn:f1_calculation}
    \begin{aligned} 
&\Tr\left(S' (K'-K)^\top (R +\gamma B^\top P_K B)(K'-K)\right) + 2 \operatorname{Tr}\left(S' (K'-K)^\top E_K \right) \\
 \leq  &4\eta_1^2 \|R +\gamma B^\top P_K B \| \operatorname{Tr}\left(
        S' E_K^\top  E_K
     \right)  - 4 \eta_1 \operatorname{Tr}\left(S' E_K^\top  E_K \right) 
   \overset{(a)}{\leq}  -2\eta_1 \operatorname{Tr}\left(S' E_K^\top  E_K\right)\\
   \leq&-2\eta_1\mu \operatorname{Tr}\left( E_K^\top  E_K\right)
  \overset{(b)}{\leq} - 2\eta_1 \mu \frac{\sigma_{\min}(R)}{\|S^{*} \|} \Tr\left( 
        S_{K^*,\Sigma^*} E_K^{\top}\left(R+\gamma B^{\top} P_K B\right)^{-1} E_K\right),
\end{aligned}\vspace{-0.5mm}
\end{equation}
where $(a)$ follows from $\eta_1 \leq  ({2 \|R +\gamma B^\top P_K B \|})^{-1}$ and $(b)$ follows from \eqref{eqn: grad K part}.

By Lemma \ref{lemma: bound of update Sigma}, $aI \preceq \Sigma' \preceq I.$ Then by Lemma \ref{lemma: smooth of fK}, 
\begin{align*}
    &f_K(\Sigma) - f_K(\Sigma')%\leq \frac{\Tr\left(  \left(R +\gamma B^\top P_K B- \frac{\tau}{2} \Sigma^{-1}\right) (\Sigma' -\Sigma)\right)}{1-\gamma} + M_{a} \Tr\left( ( \Sigma^{-1}\Sigma' - I)^2\right)
    \\
    &\leq -\frac{\eta_2}{(1-\gamma)^2} \Tr\left(\left( R +\gamma B^\top P_K B  - \frac{\tau}{2} \Sigma^{-1}\right) \Sigma \left( R +\gamma B^\top P_K B - \frac{\tau}{2} \Sigma^{-1}\right)  \Sigma  \right) \\
    &\quad \quad \quad \quad +  \frac{(\eta_2)^2 M_{a}}{(1-\gamma)^2}   \Tr\left(\left( (R +\gamma B^\top P_K B)\Sigma - \frac{\tau}{2}I\right)^2\right)\\
     \overset{(c)}{\leq}  & -\frac{\eta_2}{2(1-\gamma)^2}\Tr\left(\left( (R +\gamma B^\top P_K B)\Sigma- \frac{\tau}{2}I\right)^2\right).
\end{align*}
Here $(c)$ follows from the inequality
$
 \eta_2 =\frac{2(1-\gamma) a^2}{\tau} \leq \frac{2(1-\gamma)}{\tau} \left( \frac{\tau}{2 \| R+ \gamma B^\top P_{K} B\|}\right)^2 \overset{(d)}{\leq} \frac{2(1-\gamma)}{\tau} \overset{(e)}{\leq} \frac{1}{2M_{a}},
$
where $(d)$ is obtained from the fact that $\tau \leq 2 \sigma_{\min}(R) \leq 2 \| R+ \gamma B^\top P_{K} B\|$,  and $(e)$ follows from Lemma \ref{lemma: smooth of fK}.
 Meanwhile, observe from \eqref{eqn: f Sigma * - f Sigma} that
\begin{align*}
   f_K(\Sigma^*) - f_K(\Sigma) & \leq \frac{1}{(1-\gamma)} \Tr\Big( \left( (R+\gamma B^\top P_K B) \Sigma - \frac{\tau}{2}I\right) \cdot\\
    &\hspace{2.5cm} \Sigma^{-1} (R+\gamma B^\top P_K B)^{-1} \left( (R+\gamma B^\top P_K B) \Sigma - \frac{\tau}{2}I\right) \Big) \\
    &\leq \frac{1}{(1-\gamma)a \sigma_{\min}(R) } \Tr\left( ( (R+\gamma B^\top P_K B) \Sigma - \frac{\tau}{2}I)^2\right),
\end{align*}
while implies
\begin{equation}\label{eqn:f2_calculation}
    \begin{aligned}
    f_K(\Sigma) - f_K(\Sigma') &\leq - \frac{\eta_2 a \sigma_{\min}(R)}{2 (1-\gamma)}  \left(f_K(\Sigma) - f_K(\Sigma^*)\right).
\end{aligned}
\end{equation}
Finally, with $\zeta = \min \left\{\frac{2\mu\eta_1 \sigma_{\min}(R)}{\|S^* \|}, \frac{\eta_2 a \sigma_{\min}(R)}{2 (1-\gamma)}  \right\}$, plug 
\eqref{eqn:f1_calculation} and \eqref{eqn:f2_calculation} in \eqref{eqn: Cprime minus C} to get
$
C(K', \Sigma') - C(K,\Sigma) 
    \leq -\zeta \Big( \Tr\left( 
        S_{K^*,\Sigma^*} E_K^{\top}\left(R+\gamma B^{\top} P_K B\right)^{-1} E_K\right) + f_K(\Sigma) - f_K(\Sigma^*) \Big)
$. The proof is finished by applying \eqref{eqn: K part in C minus C *} then adding $ C(K, \Sigma) - C(K^*, \Sigma^*)$ to both sides.
%\end{proof}

\subsubsection{Proof of Lemma \ref{lemma: bounded of PK}}
%\begin{proof}{\bf of Lemma \ref{lemma: bounded of PK}}
By \eqref{eqn:C_K} and \eqref{eqn:P_K},
\begin{align*}
    C(K,\Sigma) =& \mathbb{E}_{x_0 \sim \mathcal{D}} \left[x_0^\top P_K x_0\right] + \frac{\gamma}{1-\gamma} \operatorname{Tr}(W P_K) \\
    &+ \frac{1}{1-\gamma} \operatorname{Tr}\left(\Sigma (R+\gamma B^\top P_K B)\right) -\frac{\tau}{2(1-\gamma)}(k+\log((2\pi)^k \operatorname{det}\Sigma)).
\end{align*}
Note that
\begin{align*}
    & \frac{1}{1-\gamma} \operatorname{Tr}\left(\Sigma (R+\gamma B^\top P_K B)\right) -\frac{\tau}{2(1-\gamma)}(k+\log((2\pi)^k \operatorname{det}\Sigma)) \\
    & \geq \frac{1}{1-\gamma}\left(\sigma_{\min}(R) \operatorname{Tr}(\Sigma) - \frac{\tau}{2} \left(k+k\log(2\pi)\right) 
    -\frac{\tau}{2}\log\operatorname{det}\Sigma
    \right)\\
    &\overset{(a)}{\geq}  \frac{1}{1-\gamma}\left(\frac{\tau k}{2}- \frac{\tau}{2} \left(k+k\log(2\pi)\right) 
    -\frac{\tau k}{2}\log(\frac{\tau}{2 \sigma_{\min}(R)})
    \right)= M_\tau,
\end{align*}
where $(a)$ follows from the fact that $\sigma_{\min}(R) \operatorname{Tr}(\Sigma) - \frac{\tau}{2} \left(k+k\log(2\pi)
    -\frac{\tau}{2}\log\operatorname{det}\Sigma
    \right)$ is a convex function with respect to $\Sigma$ with minimizer $\frac{\tau}{2 \sigma_{\min}(R)}I.$
Thus, 
%\begin{align*}
   $ C(K,\Sigma)  \geq  \mathbb{E}_{x_0 \sim \mathcal{D}} \left[x_0^\top P_K x_0\right] + \frac{\gamma \operatorname{Tr}(W P_K)}{1-\gamma}  + M_\tau \geq  \left(\mu + \frac{\gamma \sigma_{\min}(W)}{1-\gamma} \right) \| P_K\| + M_\tau.$
%\end{align*}
%\end{proof}

\subsection{Proofs in Section \ref{sec: policy_iteration}}
\subsubsection{Proof of Lemma \ref{lemma: contraction PI}}
% \begin{proof}{\bf of Lemma \ref{lemma: contraction PI}}
Fix $K,\Sigma$ and let $f_K$ be defined in \eqref{eqn: f}. %takes the following form: for any symmetric positive definite $X\in \mathbb{R}^{n \times n}$, \rx{do not repeat the formula if you have already defined it. Otherwise confusing.} $$ f_K(X) = \frac{\tau}{2(1-\gamma)}\log \operatorname{det} (X) -\frac{1}{1-\gamma} \operatorname{Tr} \left(X(R + \gamma B^{\top} P_K B)\right).$$
Observe that $\Sigma'$ following \eqref{alg: PI} update is the maximizer of $f_K$.
Then by Lemma \ref{lemma:almost_smooth},
    \begin{align*}
        C(&K',  \Sigma')- C(K, \Sigma) = -\operatorname{Tr} \left(S_{K', \Sigma'} E_K^\top (R+\gamma B^\top P_K B)^{-1} E_K\right) - f_K(\Sigma') + f_K(\Sigma)\\
             \overset{(a)}{\leq} & -\frac{\mu}{\| S_{K^*, \Sigma^*}\|}\Tr\left(S_{K^*, \Sigma^*} E_K^\top (R+\gamma B^\top P_K B)^{-1} E_K\right)  - \frac{\mu}{\| S_{K^*, \Sigma^*}\|} (f_K(\Sigma') - f_K(\Sigma))\\
        \overset{(b)}{\leq} &-\frac{\mu}{\| S_{K^*, \Sigma^*}\|} \Bigg( \Tr\left(S_{K^*, \Sigma^*} E_K^\top (R+\gamma B^\top P_K B)^{-1} E_K\right)  + f_K(\Sigma^*) - f_K(\Sigma)  \Bigg),
    \end{align*}
    where $(a)$ follows the fact that $f_K(\Sigma') - f_K(\Sigma) \geq 0$ and $0 <\frac{\mu}{\|S_{K^*,\Sigma^*}\|} \leq 1$; $(b)$ follows from fact that $f_K(\Sigma') \geq f_K(\Sigma^*) .$
    Finally, apply \eqref{eqn: K part in C minus C *} to get
   % \begin{align*}
       $C(K', \Sigma') - C(K, \Sigma) \leq  -\frac{\mu}{\| S_{K^*, \Sigma^*}\|}  \left(C(K,\Sigma) -  C(K^*,\Sigma^*)  \right)$
    %\end{align*}
    and add $C(K,\Sigma) -  C(K^*,\Sigma^*)$ to both sides finishes the proof.
%\end{proof}
\subsubsection{Proof of Lemma \ref{lemma: local contraction of C}}\label{sec: proof of lemma local contraction}
% \begin{proof}{\bf  of Lemma \ref{lemma: local contraction of C}}
     Let $K'$ denote $K^{(t+1)}$, $K$ denote $K^{(t)}$ and $S'$ denote $S^{(t+1)}$. Apply Lemma \ref{lemma:almost_smooth} with $f_K$ defined in \eqref{eqn: f} to get
     \begin{align*}
         &C(K', \Sigma') - C(K, \Sigma)\\
         = & - \operatorname{Tr}\left( S^*  E_K^\top (R+\gamma B^\top P_K B)^{-1} E_K\right) - \operatorname{Tr}\left(  (S' - S^*)  E_K^\top (R+\gamma B^\top P_K B)^{-1} E_K\right) \\
         & +  f_K(\Sigma) - f_K(\Sigma')\\
         \overset{(a)}{\leq} &\left(- 1 + \|(S' - S^* ) S^{* -1}\| \right) \left( \operatorname{Tr}\left( S^*  E_K^\top (R+\gamma B^\top P_K B)^{-1} E_K\right) +  f_K(\Sigma^*) - f_K(\Sigma)\right),
     \end{align*}
     where $(a)$ follows from the assumption $\| S' -S^*\| \leq \sigma_{\min}(S^*) $ which  implies $- 1 + \left\|(S' - S^* ) S^{* -1}\right\| \in [-1, 0]$ and the fact $\Sigma'$ is the maximizer of $f_K.$
            Finally, by \eqref{eqn: K part in C minus C *},
   \begin{align*}
        C(K', \Sigma') - C(K, \Sigma) &\leq \left(- 1 + \left\|(S' - S^* ) S^{* -1}\right\| \right) \left(C(K, \Sigma) - C(K^*, \Sigma^*) \right)\\
        &\leq  \left(- 1 +  \frac{\left\|S' - S^* \right\|}{\sigma_{\min}(S^*)}\right) \left(C(K, \Sigma) - C(K^*, \Sigma^*) \right).
   \end{align*}
   Adding $C(K, \Sigma) - C(K^*, \Sigma^*) $ to both sides of the above inequality finishes the proof.
\subsubsection{Proof of Lemma \ref{lemma: Sx perturbation}}\label{sec: perturbation analysis S}
This section conducts a perturbation analysis on $S_{K, \Sigma}$ and aims to prove Lemma \ref{lemma: Sx perturbation}, which bounds $\|S_{K_1, \Sigma_1} - S_{K_2, \Sigma_2}\|$ by $\left\|K_1-K_2\right\|$ and $\| \Sigma_1 - \Sigma_2\|$.

The proof of Lemma \ref{lemma: Sx perturbation} proceeds with a few technical lemmas. First, define the linear operators on symmetric matrices. For symmetric matrix $X \in \mathbb{R}^{n \times n}$, we set
\begin{equation}\label{eqn: define op}
    \begin{aligned}
    \mathcal{F}_K (X) &:= (A- BK)X(A-BK)^\top, \quad
    \mathcal{G}^K_t (X) := (A-BK)^t X (A-BK)^{\top t}\\
    \mathcal{T}_K(X) &:= \sum_{t=0}^\infty \gamma^t (A-BK)^t X (A-BK)^{\top t} = \sum_{t=0}^\infty \gamma^t \mathcal{G}^K_{t}(X).
\end{aligned}
\end{equation}
Note that when $\| A-BK\| < \frac{1}{\sqrt{\gamma}}$, we have
\begin{equation}\label{eqn: Tk as a inv}
    \mathcal{T}_K = (I-\gamma \mathcal{F}_K)^{-1}.
\end{equation}
We also define the induced norm for these operators as
$
\|T\|:=\sup _X \frac{\|T(X)\|}{\|X\|},
$
where $T=\mathcal{F}_{K}, \mathcal{G}_t^{{K}}, \mathcal{T}_{{K}}$ and the supremum is over all symmetric matrix $X \in \mathbb{R}^{n \times n}$ with non-zero spectral norm. 
We also define the induced norm for these operators as
$
\|T\|:=\sup _X \frac{\|T(X)\|}{\|X\|},
$
where $T=\mathcal{F}_{K}, \mathcal{G}_t^{{K}}, \mathcal{T}_{{K}}$ and the supremum is over all symmetric matrix $X \in \mathbb{R}^{n \times n}$ with non-zero spectral norm. 
\begin{lemma}\label{lemma: Sx}  
    $ S_{K,\Sigma} = \mathcal{T}_K (\mathbb{E}[x_0 x_0^\top]) + \sum_{t=0}^\infty \gamma^t \sum_{s=1}^t \mathcal{G}^K_{t-s} (B\Sigma B^\top + W),$ for any $ K, \Sigma \in \Omega.$
\end{lemma}
\begin{proof} (of Lemma \ref{lemma: Sx}).
Let us first show the following equation holds for any $t=1,2,3, \cdots$ by induction: %\rx{transpose missing for the last matrix?}
\begin{equation*}\label{eqn: Exx for induction}
\begin{aligned}
    \mathbb{E}[x_t x_t^\top] &= (A-BK)^t \mathbb{E}[x_0 x_0^\top] (A-BK)^{\top t} \\
    &\quad+ \sum_{s=1}^t (A-BK)^{t-s} (B\Sigma B^\top + W) (A-BK)^{\top t-s}.
 \end{aligned}
\end{equation*}
When $t=1$, it holds since 
    $$\mathbb{E}[x_1 x_1^\top] 
    = (A-BK)\mathbb{E}[x_0 x_0^\top] (A-BK)^\top + B\Sigma B^\top + W.$$
Now assume it holds for $t$, consider $t+1$:
\begin{align*}
    \mathbb{E}[x_{t+1} x_{t+1}^\top] =& (A-BK)\mathbb{E}[x_t x_t^\top](A-BK)^\top + B\Sigma B^\top + W \\
    =&  (A-BK)^{t+1} \mathbb{E}[x_0 x_0^\top] (A-BK)^{\top t+1} \\
    &  + \sum_{s=1}^t (A-BK)^{t-s +1 } (B\Sigma B^\top + W) (A-BK)^{t-s +1} + B\Sigma B^\top + W\\
    =& (A-BK)^{t+1} \mathbb{E}[x_0 x_0^\top] (A-BK)^{\top t+1}\\
    &\quad + \sum_{s=1}^{t+1} (A-BK)^{t-s +1 } (B\Sigma B^\top + W) (A-BK)^{t-s + 1}.
\end{align*}
    Thus, the induction is completed. 
    The proof is finished by noting \eqref{eqn: define op}.
    %Then by \eqref{eqn: define op}, we can write $ S_{K, \Sigma}$ as:
%$$ S_{K, \Sigma} = \sum_{t=0}^\infty\gamma^t \mathbb{E}[x_t x_t^\top] = \mathcal{T}_K\left(\mathbb{E}[x_0 x_0^\top]\right) +\sum_{t=0}^\infty \gamma^t %\sum_{s=1}^t \mathcal{G}_{t-s}^K (B\Sigma B^\top +W).$$ 
%fgdfgdfgdfgfd
\end{proof}

% {\color{red}\begin{proof} (of Lemma \ref{lemma: Sx}).
% Let us first show the following equation holds for any $t=1,2,3, \cdots$ by in\mathrm{d}uction: %\rx{transpose missing for the last matrix?}
% \begin{equation*}\label{eqn: Exx for in\mathrm{d}uction}
% \begin{aligned}
%     \mathbb{E}[x_t x_t^\top] & = (A-BK)^t \mathbb{E}[x_0 x_0^\top] (A-BK)^{\top t} \\
%     & \quad + \sum_{s=1}^t (A-BK)^{t-s} (B\Sigma B^\top + W) (A-BK)^{\top t-s}.
% \end{aligned}
% \end{equation*}
% When $t=1$, it holds since 
%     $\mathbb{E}[x_1 x_1^\top] 
%     = (A-BK)\mathbb{E}[x_0 x_0^\top] (A-BK)^\top + B\Sigma B^\top + W.$
% Now assume it holds for $t$, consider $t+1$:
% \begin{align*}
%     \mathbb{E}[x_{t+1} x_{t+1}^\top] =& (A-BK)\mathbb{E}[x_t x_t^\top](A-BK)^\top + B\Sigma B^\top + W \\
%     =&  (A-BK)^{t+1} \mathbb{E}[x_0 x_0^\top] (A-BK)^{\top t+1} \\
%     &  + \sum_{s=1}^t (A-BK)^{t-s +1 } (B\Sigma B^\top + W) (A-BK)^{t-s +1} + B\Sigma B^\top + W\\
%     =& (A-BK)^{t+1} \mathbb{E}[x_0 x_0^\top] (A-BK)^{\top t+1}\\
%     &\quad + \sum_{s=1}^{t+1} (A-BK)^{t-s +1 } (B\Sigma B^\top + W) (A-BK)^{t-s + 1}.
% \end{align*}
%     Thus, the in\mathrm{d}uction is completed. Then the proof is finished by \eqref{eqn: define op}.
% \end{proof}} 
\begin{lemma} \label{lemma: F-F'}
    $
    \left\|\mathcal{F}_{K_1}-\mathcal{F}_{K_2}\right\| \leq (\|A-BK_1\| +  \|A-BK_2\|) \|B\|\left\|K_1-K_2\right\|,
    $ for any $K_1$ and $K_2$ in $\Omega$.
\end{lemma}
Lemma \ref{lemma: F-F'} follows a direct calculation with the triangle inequality.
%Note that Lemma \ref{lemma: F-F'} follows a direct calculation with triangle inequality.
\begin{lemma}\label{lemma: G-G'}

For any $K_1 \in \Omega$ and $K_2 \in \Omega$, with $\xi_{\gamma, \rho}$ defined in \eqref{eqn: constants xi zeta omega}, 
   \begin{align}
   \left\|\left(\mathcal{T}_{K_1} - \mathcal{T}_{K_2}  \right)(X)\right\| &\leq \sum_{t=0}^\infty \gamma^t \left\| \left(\mathcal{G}_t^{K_1} - \mathcal{G}^{K_2}_t \right) (X)\right\| \leq  \xi_{\gamma, \rho}\|\mathcal{F}_{K_1} - \mathcal{F}_{K_2} \| \|X\|, \label{eqn: T diff}\\
         \sum_{t= 0}^{T-1} \| (\mathcal{G}_{t}' -\mathcal{G}_{t})(X) \| 
&\leq \frac{2-\rho^2 -  \rho^{2T}}{(1-\rho^2)^2 } \|\mathcal{F'} -\mathcal{F} \| \| X\|, \quad \forall \, T \geq 1. \label{eqn: G diff}
    \end{align}

\end{lemma}
%{Proofs of Lemma \ref{lemma: Sx}, \ref{lemma: F-F'}, and \ref{lemma: G-G'} can be found in the online companion.}
\begin{proof}{ (of Lemma \ref{lemma: G-G'}).}
  To ease the exposition,  let  us denote $\mathcal{G}_t = \mathcal{G}^{{K_1}}_t, \mathcal{G}'_t = \mathcal{G}^{{K_2}}_t, \mathcal{F} = \mathcal{F}_{K_1}$ and $\mathcal{F}' = \mathcal{F}_{K_2}$. Then for any symmetric matrix $X \in \mathbb{R}^{n \times n}$ and $t = 0,1,2, \cdots$,
    \begin{equation}\label{eqn: lemma G}
        \begin{aligned}
\left\|\left(\mathcal{G}_{t+1}^{\prime}-\mathcal{G}_{t+1}\right)(X)\right\| 
% & =\left\|\mathcal{F}^{\prime} \circ \mathcal{G}_t^{\prime}(X)-\mathcal{F} \circ \mathcal{G}_t(X)\right\| \\
& =\left\|\mathcal{F}^{\prime} \circ \mathcal{G}_t^{\prime}(\Sigma)-\mathcal{F}^{\prime} \circ \mathcal{G}_t(X)+\mathcal{F}^{\prime} \circ \mathcal{G}_t(X)-\mathcal{F} \circ \mathcal{G}_t(X)\right\| \\
& \leq\left\|\mathcal{F}^{\prime} \circ \mathcal{G}_t^{\prime}(X)-\mathcal{F}^{\prime} \circ \mathcal{G}_t(X)\right\|+\left\|\mathcal{F}^{\prime} \circ \mathcal{G}_t(X)-\mathcal{F} \circ \mathcal{G}_t(X)\right\| \\
& =\left\|\mathcal{F}^{\prime} \circ\left(\mathcal{G}_t^{\prime}-\mathcal{G}_t\right)(X)\right\|+\left\|\left(\mathcal{F}^{\prime}-\mathcal{F}\right) \circ \mathcal{G}_t(X)\right\| \\
& \leq\left\|\mathcal{F}^{\prime}\right\|\left\|\left(\mathcal{G}_t^{\prime}-\mathcal{G}_t\right)(X)\right\|+\left\|\mathcal{G}_t\right\|\left\|\mathcal{F}^{\prime}-\mathcal{F}\right\|\|X\| \\
& \leq \rho^2\left\|\left(\mathcal{G}_t^{\prime}-\mathcal{G}_t\right)(X)\right\|+\rho^{2t}\left\|\mathcal{F}^{\prime}-\mathcal{F}\right\|\|X\|.
\end{aligned}
    \end{equation}
    Summing \eqref{eqn: lemma G} for $t=0, 1, 2, \cdots$ with a discount factor $\gamma$, we have %\rx{inequality below?}
    \begin{align*}
      \sum_{t=0}^\infty \gamma^t \left\|\left(\mathcal{G}_{t+1}^{\prime}-\mathcal{G}_{t+1}\right)(X)\right\| \leq  \rho^2 \sum_{t=0}^\infty \gamma^t \left\|\left(\mathcal{G}_{t}^{\prime}-\mathcal{G}_{t}\right)(X)\right\| + \frac{1}{1- \gamma \rho^2}\left\|\mathcal{F}^{\prime}-\mathcal{F}\right\|\|X\|.
    \end{align*}
   Note that the left-hand side of the above inequality equals to
    $$
        \frac{1}{\gamma} \sum_{t=0}^\infty \gamma^{t+1}\left\|\left(\mathcal{G}_{t+1}^{\prime}-\mathcal{G}_{t+1}\right)(X)\right\|  =   \frac{1}{\gamma} \sum_{t=0}^\infty \gamma^{t}\left\|\left(\mathcal{G}_{t}^{\prime}-\mathcal{G}_{t}\right)(X)\right\| - \frac{1}{\gamma} \left\|\left(\mathcal{G}_{0}^{\prime}-\mathcal{G}_{0}\right)(X)\right\|.
    $$
    Thus, by combining the above calculations, we obtain:
    \begin{align*}
          \left(  \frac{1}{\gamma} -\rho^2 \right)\sum_{t=0}^\infty \gamma^{t}\left\|\left(\mathcal{G}_{t}^{\prime}-\mathcal{G}_{t}\right)(X)\right\| &\leq \frac{1}{\gamma} \left\|\left(\mathcal{G}_{0}^{\prime}-\mathcal{G}_{0}\right)(X)\right\| + \frac{1}{1- \gamma \rho^2}\left\|\mathcal{F}^{\prime}-\mathcal{F}\right\|\|X\|. 
    \end{align*}
    The proof of \eqref{eqn: T diff} is finished by noting $\mathcal{G}_0 = \mathcal{F}$ and $\mathcal{G}'_0 = \mathcal{F}'$. 

To show \eqref{eqn: G diff}, fix integer $T \geq 1.$ Summing up \eqref{eqn: lemma G} for $t = 0, 1, \cdots, T-1$ to get
 \begin{align}\label{eqn: G diff finite}
      \sum_{t=0}^{T-1} \left\|\left(\mathcal{G}_{t+1}^{\prime}-\mathcal{G}_{t+1}\right)(X)\right\| \leq  \rho^2 \sum_{t=0}^{T - 1}  \left\|\left(\mathcal{G}_{t}^{\prime}-\mathcal{G}_{t}\right)(X)\right\| + \frac{1-\rho^{2 T}}{1-\rho^2}\left\|\mathcal{F}^{\prime}-\mathcal{F}\right\|\|X\|.
    \end{align}
    Note that the left-hand side equals 
$$
        \sum_{t=0}^{T-1} \left\|\left(\mathcal{G}_{t+1}^{\prime}-\mathcal{G}_{t+1}\right)(X)\right\|  =  \sum_{t=0}^{T-1}\left\|\left(\mathcal{G}_{t}^{\prime}-\mathcal{G}_{t}\right)(X)\right\| - \left\|\left(\mathcal{F}^{\prime}-\mathcal{F}\right)(X)\right\|.
    $$
Plugging the above equation into \eqref{eqn: G diff finite} finishes the proof of \eqref{eqn: G diff}.
\end{proof}

\begin{proof} (of Lemma \ref{lemma: Sx perturbation}). Denote $\mathcal{G}_t = \mathcal{G}^{{K_1}}_t, \mathcal{G}'_t = \mathcal{G}^{{K_2}}_t, \mathcal{T} = \mathcal{T}_{K_1}, \mathcal{T}' = \mathcal{T}_{K_2}, \mathcal{F} = \mathcal{F}_{K_1}$ and $\mathcal{F}' = \mathcal{F}_{K_2}$ to ease the exposition.
Observe that \vspace{-2mm}
\begin{equation*}
    \begin{aligned}
 & \|S_{K_1, \Sigma_1} - S_{K_2, \Sigma_2}\| \overset{(a)}{\leq} \| (\mathcal{T} - \mathcal{T}' )\mathbb{E}[x_0 x_0^\top] \| +\| \sum_{t=0}^\infty \gamma^t \sum_{s=1}^t (\mathcal{G}_{s}^{'} -\mathcal{G}_{s} )(B\Sigma_1 B^\top +W) \| \\
    &\qquad\qquad\qquad\qquad\qquad +  \sum_{t=0}^\infty  \gamma^t \sum_{s=1}^t \|\mathcal{G}_{t-s}^{\prime}  (B\Sigma_1 B^\top + W)  -  \mathcal{G}_{t-s}^{\prime}(B\Sigma_2 B^\top + W)\| \\
% &\overset{(a)}{\leq} \|  (\mathcal{T} - \mathcal{T}' )(\mathbb{E}[x_0 x_0^\top]) \|  + \| \sum_{t=0}^\infty \gamma^t \sum_{s=0}^{t-1} (\mathcal{G}_{s}^{'} -\mathcal{G}_{s} )(B\Sigma_1 B^\top +W) \| \\
% &\quad\quad+ \sum_{t=0}^\infty \gamma^t \sum_{s=1}^t \rho^{2(t-s)} \| B\|^2 \|\Sigma_1 - \Sigma_2 \|\\
&\overset{(b)}{\leq}\xi_{\gamma, \rho}  \|\mathcal{F} - \mathcal{F}' \| \|\mathbb{E}[x_0 x_0^\top]\| 
+ \sum_{t=0}^\infty\gamma^t \frac{2-\rho^2 - \rho^{2t}}{(1-\rho^2)^2}\|\mathcal{F'} -\mathcal{F} \| \| B\Sigma_1 B^\top +W\| \\
&\qquad+ 
    \sum_{t=0}^\infty \gamma^t \sum_{s=1}^t \rho^{2(t-s)} \| B\|^2 \| \Sigma_1 - \Sigma_2\|,
    % \\
    % &= \xi_{\gamma, \rho} \|\mathcal{F} - \mathcal{F}{'} \| \|\mathbb{E}[x_0 x_0^\top]\| +\zeta_{\gamma, \rho} \|\mathcal{F'} -\mathcal{F} \| \cdot \|B\Sigma_1 B^\top +W\|+ \omega_{\gamma, \rho}\| B\|^2 \| \Sigma_1 - \Sigma_2\| ,
\end{aligned}\vspace{-2mm}
\end{equation*}
where $(a)$ is from \eqref{eqn: def S K Sigma} and $(b)$ follows from  Lemma \ref{lemma: G-G'} and 
$\| \mathcal{G}_t(X) -  \mathcal{G}_t(X') \| = \| (A-BK)^t (X-X') (A-BK)^{\top t}\| \leq \rho^{2t} \| X - X'\|,$
$\forall X, X' \in\mathbb{R}^{n \times n}, \forall t \geq 0 .$
Finally, applying \eqref{eqn: constants xi zeta omega} and Lemma \ref{lemma: F-F'} finishes the proof.
\end{proof}

\subsubsection{Proof of Lemma \ref{lemma: Sprime minus Sstar}} 
The objective of this section is to provide a bound for $\|S^{(t+1)} - S^*\|$ in terms of $\| K^{(t)} - K^*\|$, as summarized in Lemma \ref{lemma: Sprime minus Sstar}.

Note that
Lemma \ref{lemma: Sx perturbation} can be employed to derive a bound on $\|S^{(t+1)} - S^* \|$ in relation to $\| K^{(t+1)} - K^*\|$ and $\| \Sigma^{(t+1)} - \Sigma^*\|.$ In this section, we further establish this bound by deriving bounds for $\| K^{(t+1)} - K^*\|$ and $\| \Sigma^{(t+1)} - \Sigma^*\|$ in terms of $\| K^{(t)} - K^*\|$  and $\| P_{K^{(t)}} - P_{K^*}\|$ (\textit{cf.} Lemma \ref{lemma: K' - K} and Lemma \ref{lemma: Sigma'-Sigma_star}). Additionally, the perturbation analysis for $P_K$ in Lemma \ref{lemma: P minus Pprime} demonstrates $\| P_{K^{(t)}} - P_{K^*}\|$ can be bounded by $\| {K^{(t)}} - {K^*}\|$, which completes the proof for Lemma \ref{lemma: Sprime minus Sstar}.
 % We first establish the bound for the one-step update of $K$ and $\Sigma$ in Lemma \ref{lemma: K' - K} and Lemma \ref{lemma: Sigma'-Sigma_star} respectively.
\begin{lemma}[Bound of one-step update of $K$]\label{lemma: K' - K}
Assume the update of parameter $K$ follows the updating rule in \eqref{alg: PI}. Then it holds that: 
\begin{equation*}
    \begin{aligned}
        \| K^{(t+1)}-K^{*}\| \leq & \left(1+ \sigma_{\min}(R)  \|\gamma B^{\top} P_{K^*}B + R \| \right) \cdot  \|K^{(t)}- K^* \|  \\
       & + \gamma \sigma_{\min}(R) \left( \| B\| \| A\| + \| B\|^2 \kappa \right) \cdot \|P_{K^{(t)}} - P_{K^*} \|.
    \end{aligned}\vspace{-2mm}
\end{equation*}
\end{lemma}

\begin{proof} (of Lemma \ref{lemma: K' - K}). Let $K'$ denote $K^{(t+1)}$ and $K$ denote $K^{(t)}$ to ease the notation.
Theorem \ref{thm: optimal val func and policy} shows that for an optimal $K^*$,
$$
K^* = \gamma(R+\gamma B^{\top} P_{K^*} B)^{-1} B^{\top} P_{K^*} A.
$$
Then, by the definition of $E_K$ in Lemma \ref{thm: PG},
\begin{equation}%\label{eqn: E_K *}
    E_{K^*} =-\gamma B^{\top} P_{K^*}A+ (\gamma B^{\top} P_{K^*}B + R )K^* = 0.
\end{equation}
Now we bound the difference between $K' - K$: 
    \begin{eqnarray}
         \|K' -K^* \| &=& \| K - (R+\gamma B^\top P_K B)^{-1} E_K - K^* + (R+\gamma B^\top P_{K} B)^{-1} E_{K^*}\| \nonumber\\
         &\leq& \|K- K^* \|%+ \|  (R+\gamma B^\top P_K B)^{-1} E_K - (R+\gamma B^\top P_{K} B)^{-1} E_{K^*}\|\\
          + \| (R+\gamma B^\top P_{K^*} B)^{-1} \| \|  E_{K^*} - E_{K^*}\|\nonumber
         \\
         &\leq& \|K- K^* \| + \sigma_{\min}(R)\|E_K - E_{K^*} \| .\label{eq:int1}
    \end{eqnarray}
    To bound the difference between $E_K$ and $E_{K^*}$, observe: 
    \begin{eqnarray}
        \|E_K - E_{K^*} \| &\leq& \gamma \| B\| \| A\| \| P_{K} - P_{K^*}\| + \| (\gamma B^{\top} P_{K^*}B + R )K^* - (\gamma B^{\top} P_{K^*}B + R )K\| \nonumber\\
        &&+ \| (\gamma B^{\top} P_{K^*}B + R )K - (\gamma B^{\top} P_{K}B + R )K\|\nonumber\\
        &\leq &  \gamma \| B\| \| A\| \| P_{K} - P_{K^*}\| + \|\gamma B^{\top} P_{K^*}B + R \| \|K^* - K\|\nonumber\\
        &&+ \gamma \| B\|^2 \| P_{K^*}-P_{K}\| \|K\|.\label{eq:int2}
    \end{eqnarray}
    Combining \eqref{eq:int1} and \eqref{eq:int2}, then using $\kappa \geq \| K\|$ for any $K\in\Omega$ completes the proof.
\end{proof}
\begin{lemma}[Bound of one-step update of $\Sigma$]\label{lemma: Sigma'-Sigma_star}
Suppose $\{ K^{(t)}, \Sigma^{(t)} \}_{t=0}^\infty$ follows the update rule in \eqref{alg: PI}. Then we have
    $
    \| \Sigma^{(t+1)} -  \Sigma^{*}\| \leq \frac{\tau \gamma \|B\|^2}{\sigma_{\min}(R)} \| P_{K^{(t)}} - P_{K^*}\|.
    $
    
   \begin{proof} (of Lemma \ref{lemma: Sigma'-Sigma_star}). Observe from \eqref{alg: PI} and \eqref{eqn: Sigma *} that 
    \begin{align*}
       & \left\| \Sigma^{(t+1)} - \Sigma^* \right\| = \frac{\tau}{2} \left\|(R+\gamma B^\top P_{K^{(t)}} B)^{-1} - (R+\gamma B^\top P_{K^*} B)^{-1} \right\|\\
        &=\frac{\tau}{2}  \left\|(R+\gamma B^\top P_{K^{(t)}} B)^{-1} \cdot \gamma B^\top \left( P_{K^{(t)}} - P_{K^*}\right) B \cdot (R+\gamma B^\top P_{K^*} B)^{-1}\right\|\\
        &\leq \frac{\tau \gamma \|B\|^2}{ 2  \sigma_{\min}(R)^2} \left\| P_{K^{(t)}} - P_{K^*}\right\|.
    \end{align*}
\end{proof}

\end{lemma}
\begin{lemma}[$P_K$ perturbation]\label{lemma: P minus Pprime}For any $K \in \Omega,$ with $c$ defined in \eqref{eqn: define C1 C2},
    $  \|P_{K}-P_{K^*} \| \leq c \| K - K^*\|.$ %{\color{red} proofs online}
\end{lemma}

\begin{proof} (of Lemma \ref{lemma: P minus Pprime}).
Fix $K \in \Omega.$ 
    % Fix $K_1$ and $K_2$ such that $\|A-BK_1\| < \frac{1}{\sqrt{\gamma}}$ and  $\|A-BK_2\| < \frac{1}{\sqrt{\gamma}}$. 
    By \eqref{eqn: Tk as a inv} and \eqref{eqn:P_K}, %we have for any $K$, 
    %\begin{align*}
        $$\mathcal{T}_{K}^{-1}(P_K) = (I-\gamma \mathcal{F}_K) (P_K)= P_K - \gamma (A-BK)^\top P_K (A-BK)= Q+ K^\top R K,$$
    %\end{align*}
    which immediately implies $P_K = \mathcal{T}_K (Q+K^\top R K).$ Similarly $P_{K^*} = \mathcal{T}_{K^*} (Q+{K^*}^\top R {K^*}).$
    To bound the difference between $P_{K}$ and $P_{K^*}$, observe that,
    \begin{align}
    & \|P_{K}-P_{K^*} \| 
    % =  \left\|\mathcal{T}_{{K_2}}\left(Q+ K_2^{\top} R {K_2} \right)-\mathcal{T}_{K_1}\left(Q+{K_1}^{\top} R {K_1}\right)\right\| \notag \\ 
    % \leq & \Big\| \mathcal{T}_{K_2}\left(Q+K_2^{\top} R K_2\right)-\mathcal{T}_{K_1} \left(Q+\left(K_2\right)^{\top} R K_2\right)  \notag\\ 
    % & -\left(\mathcal{T}_{K_1}\left(Q+{K_1}^{\top} R {K_1}\right)-\mathcal{T}_{K_1}\left(Q+K_2^{\top} R K_2\right)\right) \Big\|  \notag \\ 
    % = & \left\|\mathcal{T}_{K_2}\left(Q+K_2^{\top} R K_2\right)-\mathcal{T}_{K_1}\left(Q+K_2^{\top} R K_2\right)-\mathcal{T}_{K_1} \circ\left({K_1}^{\top} R {K_1}-K_2^{\top} R K_2\right)\right\|  \notag\\
    \leq  \left\|\mathcal{T}_{K}\left(Q+K^{\top} R K\right)-\mathcal{T}_{K^*}\left(Q+K^{\top} R K\right)\right\| 
    + \| \mathcal{T}_{K^*}\|\| {K^*}^{\top} R {K^*}-K^{\top} R K\|.  \label{eqn: part 0 in P- P'}
    \end{align}
    For the first term in \eqref{eqn: part 0 in P- P'}, we can apply Lemma \ref{lemma: G-G'} and Lemma \ref{lemma: F-F'} to get
   \begin{equation}\label{eqn: part 1 in P- P'}
        \begin{aligned}
&\left\|\mathcal{T}_{K}\left(Q+K^{\top} R K\right)-\mathcal{T}_{K^*}\left(Q+K^{\top} R K\right)\right\|  \\
        \leq &\xi_{\gamma \rho} \|\mathcal{F}_{K^*} - \mathcal{F}_{K} \| \|Q+K^{\top} R K \| 
        \leq 2 \rho \xi_{\gamma \rho} \|B\|\left\|K^*-K\right\| \cdot \left(\|Q\| + \| R\| \|K\|^2 \right).
        \end{aligned}
    \end{equation}
For the second term in \eqref{eqn: part 0 in P- P'}, 
note that by Lemma 17 in \cite{fazel2018global}, 
    $
    \| \mathcal{T}_K\| \leq \frac{1}{\mu} \| \mathcal{T}_K(\mathbb{E}[x_0 x_0^\top])\|.
    $
    Since $S_{K,\Sigma} \succeq \mathcal{T}_K(\mathbb{E}[x_0 x_0^\top])$, thus 
    $   \| \mathcal{T}_K\| \leq 
    % \frac{1}{\mu} \| \mathcal{T}_K(\mathbb{E}[x_0 x_0^\top])\|=
    \frac{1}{\mu} \sigma_{max}\left( \mathcal{T}_K(\mathbb{E}[x_0 x_0^\top])\right) \leq \frac{1}{\mu} \| S_{K,\Sigma} \|. $
Then we have
\begin{equation}\label{eqn: part 2 in P- P'}
    \begin{aligned}
        & \| \mathcal{T}_{K^*}\|\| {K^*}^{\top} R {K^*}-K^{\top} R K \| \\
   =&   \| \mathcal{T}_{K^*}\|\| {K^*}^{\top} R {K^*}- {K^*}^\top R K + {K^*}^\top R K- K^{ \top} R K \| \\
   \leq& \| \mathcal{T}_{K^*}\| \| R \| \| K-{K^*}\| \left( \| {K^*}\|  + \| K\|\right)\\
   \leq &  \frac{\| S_{{K^*},\Sigma^*} \| }{\mu}  \| R \| \| K-{K^*}\|\left( \| {K^*}\|  + \| K\|\right). 
    \end{aligned}
\end{equation}
Plugging \eqref{eqn: part 1 in P- P'}, \eqref{eqn: part 2 in P- P'}, and \eqref{eqn: define C1 C2} in \eqref{eqn: part 0 in P- P'} completes the proof.
\end{proof}
With these lemmas, the proof of Lemma \ref{lemma: Sprime minus Sstar} is completed as follows:
\begin{proof} (of Lemma \ref{lemma: Sprime minus Sstar}). Let $K'$ denote $K^{(t+1)}$ and $K$ denote $K^{(t)}$.
% Apply Lemma \ref{lemma: P minus Pprime} with $K$ and $K^*$ and \eqref{eqn: define C1 C2} to get
%  \begin{align}
%         & \|P_{K^*}-P_{K} \| \leq c \left\|K -K^*\right\|.\label{eqn: P* - PK in a lemma} 
%     \end{align}
  With the assumption that $\|A-BK' \| \leq \rho$,  $\|A-BK^*\| \leq \rho$, we can apply Lemma \ref{lemma: Sx perturbation} to get\\
  \begin{equation}\label{eqn: S* - S' in a lemma}
       \begin{aligned}
      \|&S_{K^*, \Sigma^*}   - S_{K', \Sigma'}\| 
      \leq   \omega_{\gamma, \rho}\| B\|^2 \| \Sigma^* - \Sigma'\|  \\
     &\quad + \left( \xi_{\gamma, \rho}  \cdot 
     \| \mathbb{E}[x_0 x_0^\top]\|+ \zeta_{\gamma, \rho}  \cdot \|B\Sigma^* B^\top +W\| \right) \cdot 2 \rho \|B\|\left\|K^*-K'\right\|.
\end{aligned}
  \end{equation} 
Apply Lemma \ref{lemma: K' - K} and Lemma \ref{lemma: P minus Pprime} to get
\begin{equation*}\label{eqn: K'-K* in a lemma}
    \begin{aligned}
          \| K'-K^*\| 
    \leq   \Big(&1+ \sigma_{\min}(R) \cdot \|\gamma B^{\top} P_{K^*}B + R \| + c \gamma \sigma_{\min}(R) \left( \| B\| \| A\| + \| B\|^2 \kappa  \right) \Big) \\
    &\cdot \|K- K^* \|.  
    \end{aligned}
\end{equation*}
Apply Lemma \ref{lemma: Sigma'-Sigma_star} and \ref{lemma: P minus Pprime} to get
$    \| \Sigma^{\prime} -  \Sigma^{*}\| \leq \frac{\tau \gamma \|B\|^2  \| P_{K} - P_{K^*}\|}{2 \sigma_{\min}(R)^2} \leq c \frac{\tau \gamma \|B\|^2 \| {K} - {K^*}\|}{2 \sigma_{\min}(R)^2}  .$
Finally, plugging the above two inequalities
into \eqref{eqn: S* - S' in a lemma} finishes the proof.
\end{proof}

\paragraph{Acknowledgement} Renyuan Xu is partially supported by the NSF CAREER Award DMS-2524465.
\bibliographystyle{siamplain}
\bibliography{ref}

\end{document}